\RequirePackage{ifpdf}
\ifpdf 
\documentclass[pdftex]{sigma}
\else
\documentclass{sigma}
\fi

\numberwithin{equation}{section}
\numberwithin{theorem}{section}
\numberwithin{proposition}{section}
\numberwithin{lemma}{section}
\numberwithin{corollary}{section}
\numberwithin{definition}{section}
\numberwithin{example}{section}
\numberwithin{remark}{section}
\numberwithin{note}{section}

\begin{document}


\renewcommand{\thefootnote}{$\star$}

\renewcommand{\PaperNumber}{103}

\FirstPageHeading

\ShortArticleName{Isomorphism of Intransitive Linear Lie Equations}

\ArticleName{Isomorphism of Intransitive Linear Lie Equations\footnote{This paper is a
contribution to the Special Issue ``\'Elie Cartan and Dif\/ferential Geometry''. The
full collection is available at
\textit{}\href{http://www.emis.de/journals/SIGMA/Cartan.html}{http://www.emis.de/journals/SIGMA/Cartan.html}}}

\Author{Jose Miguel Martins VELOSO}

\AuthorNameForHeading{J.M.M. Veloso}

\Address{Faculdade de Matematica, UFPA, Belem, PA, CEP 66075-110, Brasil}
\Email{\href{mailto:veloso@ufpa.br}{veloso@ufpa.br}}

\ArticleDates{Received February 09, 2009, in f\/inal form November 11, 2009;  Published online November 17, 2009}

\Abstract{We show that formal isomorphism of intransitive linear Lie equations along transversal to the orbits can be extended to neighborhoods of these transversal. In analytic cases, the word formal is dropped from theorems. Also, we associate an intransitive Lie algebra with each intransitive linear Lie equation, and from the intransitive Lie algebra we recover the linear Lie equation, unless of formal isomorphism. The intransitive Lie algebra gives the structure functions introduced by \'E.~Cartan.}

\Keywords{Lie equations; Lie groupoids; intransitive; isomorphism}

\Classification{58H05; 58H10}

\newcommand{\zerobra}[2]{{{[}#1,#2{]}}}
\newcommand{\firstbra}[2]{{{[\![}#1,#2{]\!]}}}
\newcommand{\secondbra}[2]{{\rule[-3.03pt]{2pt}{12.05pt}{\hskip-2pt}{\bm{[}}{#1},{#2}{\bm{]}{\hskip-2pt}\rule[-3.03pt]{2pt}{12.05pt}}}}
\newcommand{\thirdbra}[2]{{\rule[-3.03pt]{3.5pt}{12.05pt}{\hskip-2pt}{\bm{[}}{#1},{#2}{\bm{]}{\hskip-2pt}\rule[-3.03pt]{3.5pt}{12.05pt}}}}

\newcommand{\bm}[1]{\mbox{\boldmath $#1$}}

\section{Introduction}

It is known \cite{G4} that the isomorphism of the f\/ibers of transitive linear Lie equations at two points is suf\/f\/icient to obtain the formal isomorphism of the Lie equations. This is proved by constructing a system of partial dif\/ferential equations (SPDE) whose solutions would be these isomorphisms. This SPDE may  be not integrable \cite{Co,GuSt2}, although it is formally integrable. If the data are analytic, the SPDE is integrable.

In this paper we consider principally the extension of this theorem to intransitive linear Lie equations. Intransitive linear Lie equations generate a family of orbits on the manifold  that we suppose locally to be a foliation.
Given two intransitive Lie equations, consider the restriction of both Lie equations along  two transversal to the orbits. If these restrictions are isomorphic in a certain sense, then we can construct a formally integrable SPDE whose solutions (if they exist) are isomorphisms of the two Lie equations. Therefore we prove, at least in the analytic case, that formal isomorphism of two linear Lie equations along transversal to the orbits can be extended locally to local isomorphism of the two linear Lie equations in a neighborhood of the two transversal.

Specif\/ically, consider $M$ and $M'$  manifolds of the same dimension, $V$ and $V'$ integrable distributions of the same dimension on $M$ and $M'$, respectively, $N$ and $N'$ submanifolds of $M$ and $M'$ such that each integral leaf of $V$ and $V'$ through points $x\in M$ and $x'\in M'$ intersect $N$ and $N'$ (at least locally) at unique points $\rho x$ and $\rho'x'$, respectively. Let be $\phi:N\rightarrow N'$  a local dif\/feomorphism, $a\in N$, $a'\in N'$, $\phi(a)=a'$, $Q^k_\phi$ the manifold of $k$-jets of local dif\/feomorphisms $f:M\rightarrow M'$ such that $\phi\rho(x)=\rho'f(x)$, and $\mathcal{Q}^k_\phi$ the sheaf of germs of invertible local sections of $Q^k_{\phi}$. Furthermore, let be $R^k\subset J^kV$, $R'^k \subset J^kV'$ intransitive linear Lie equations such that $R^0=J^0V$, $R'^0=J^0V'$. We say that $R^k$ at point $a$ is \emph{formally isomorphic} to $R'^k$ at a point $a'$ if we can construct a formally integrable SPDE $S^k\subset Q^k_\phi$ such that any solution $f$ of $S^k$ satisf\/ies $(j^{k+1}f)_*(R^k)=R'^k$ (see Def\/inition~\ref{foriso}). We prove in Proposition~\ref{FQ} that this condition is equivalent to the existence of  $F\in \mathcal{Q}^{k+1}_\phi$ such that $\beta F(a)=a'$ and $F_*(TM\oplus R^k)=TM'\oplus R'^k$. Then we prove:\vspace{.3cm}

\noindent \textbf{Theorem~\ref{iso}.}
\textit{Let be $\Phi:N\rightarrow Q^{k+1}_\phi$ such that $\beta\Phi=\phi$ and
\[\Phi_*\big(TN\oplus R^k|_N\big)=TN'\oplus R'^k|_{N'}.\]
Then given a diffeomorphism $f:M\rightarrow M'$ such that $f_*V=V'$ and $f|_N=\phi$ there exists $F\in\mathcal Q^{k+1}_\phi$ satisfying $F|_N=\Phi$, $\beta F=f$ and}
\[F_*\big(TM\oplus R^k\big)=TM'\oplus R'^k.\]
 This theorem for transitive linear Lie equations is in  \cite{G4}.

The next step is to def\/ine an intransitive Lie algebra representing $TN\oplus R^k|_N$. Transitive Lie algebras were def\/ined in \cite{GuSt1} as the algebraic object necessary to study transitive inf\/initesimal Lie pseudogroups. This program was pursued in the papers \cite{G3,G4,G5,GS,KS,Ma1,Ma2,RP,SS} disclosing the f\/itness of transitive Lie algebras to study  transitive linear Lie equations.

At the same time, several tries were made to include intransitive linear Lie equations in this theory \cite{Ki1,Ki2,Mo,N}.  Basically, they associated a~family of transitive Lie algebras on a~transversal to the orbits with an intransitive inf\/initesimal Lie pseudogroup, each transitive Lie algebra corresponding to the transitive inf\/initesimal Lie pseudogroup obtained  by restriction to each orbit of the inf\/initesimal intransitive Lie pseudogroup. Another approach was to study the Lie algebra of inf\/inite jets at one point of the inf\/initesimal intransitive Lie pseudogroup; this algebra is bigraded, and this bigraduation may give the power series of the structure functions introduced in \cite{C}. In the f\/irst approach, it may happen to be impossible to relate transitive Lie algebras along the transversal and give a continuous structure to this family, as the following example in \cite{V3} shows: let be $\Theta$ the inf\/initesimal intransitive Lie pseudogroup  acting in the plane~$\mathbb{R}^2$ given by
\[\Theta=\left\{\theta(x,y)\frac{\partial}{\partial y}\mbox{ such that }\frac{\partial\theta}{\partial y}=a(x)\frac{\partial\theta}{\partial x}\mbox{ with }a\neq 0, \ a(0)=0\right\}.\]
Two of these inf\/initesimal intransitive Lie pseudogroups given by functions $a(x)$ and $\bar a(x)$ are isomorphic if and only if there exists a $C^\infty$-function $b(x)$, with $b(0)\neq 0$ such that $\bar a(x)=b(x)a(x)$. The restriction of $\Theta$ to the orbits $\{(x_0,y):y\in\mathbb{R}\}$, with $a(x_0)\neq 0$ is the inf\/initesimal Lie pseudogroup of dif\/ferentiable vector f\/ields on $\mathbb{R}$, and the restriction to the orbit where \mbox{$a(x_0)=0$} are the inf\/initesimal translations on $\mathbb{R}$. This example shows that, even if the intransitive linear Lie equation associated to $\Theta$ has all properties of regularity, by choosing an appropriate function~$a(x)$, we obtain a highly discontinuous family of transitive Lie algebras along the transversal. If $a(x)$ and $\bar a(x)$ have the same power expansion series around $0$, the bigraded Lie  algebras of inf\/inite jets of vector f\/ields of $\Theta$ and $\bar\Theta$ at point $(0,0)$ are the same. So, the bigraded Lie algebras cannot distinguish between  non isomorphic intransitive Lie pseudogroups. However, \'E.~Cartan associated ``structure functions" constant along the orbits with inf\/initesimal Lie pseudogroups. Therefore any def\/inition of intransitive Lie algebras must contain a way of getting the germs of ``structure functions" at the point considered.

A way to def\/ine such algebra is the following: let be $\mathcal{O}_{N,a}$ the ring of germs of real functions def\/ined on $N$ at the point $a$,  $L_m$  the $\mathcal{O}_{N,a}$-module of germs of sections of $R^m|_N$ at point~$a$, and~$\mathcal TN_a$ the $\mathcal{O}_{N,a}$-module of germs of sections of $TN$ at point $a$. The vector bundles  $R^m$, with $m> k$, are the prolongations of $R^k$, which we suppose formally integrable. There is a $\mathbb{R}$-bilinear map
\[\firstbra{\;}{\;}_m: \ (\mathcal{T}N_a\oplus L_m)\times (\mathcal{T}N_a\oplus L_m)\rightarrow (\mathcal{T}N_a\oplus L_{m-1}),\]
given by: the bracket of germs of vector f\/ields in $\mathcal{T}N_a$; the germ of $i(v)D\xi$ for $v\in\mathcal{T}N_a$, $\xi \in L_m$ and $D$ the linear Spencer operator; and the $\mathcal{O}_{N,a}$-bilinear map def\/ined on holonomic sections by $\firstbra{j^k\theta}{j^k\eta}_k=j^{k-1}[\theta,\eta]$. Let be $L=\mbox{limproj}_mL_m$, $\mathcal{L}=\mathcal{T}N_a\oplus L$ and $\firstbra{\;}{\;}_\infty=\mbox{limproj}_m\firstbra{\;}{\;}_m$. Then $(\mathcal{L},\firstbra{\;}{\;}_\infty)$ is the $\mathcal{O}_{N,a}$-\emph{intransitive} $\mathbb{R}$-\emph{Lie algebra} associated with the formally integrable linear Lie equation $R^k$ at the point $a$. From  $(\mathcal{L},\firstbra{\;}{\;}_\infty)$ we can obtain the germ of $R^k|_N$ at point $a$, and, by applying the Theorem of \cite{V} and Theorem~\ref{iso} of this paper, we get the germ of the linear Lie equation $R^k$ at the point $a$, up to a formal isomorphism (cf. Theorems~\ref{VV} and \ref{isore} below).

We summarize the content of this paper. Section~\ref{pre} presents basic facts on groupoids and algebroids of jets, the calculus on the diagonal introduced in \cite{Ma1,Ma2}, and the construction of the f\/irst linear and non-linear Spencer complexes. We tried to be as complete as possible, and the presentation emphasizes the geometric relationship  between the f\/irst linear and non-linear Spencer operators, and the left and right actions of a groupoid on itself. We hope this section will facilitate the reading of this paper, since several formulas given here are sometimes not easily identif\/ied in \cite{Ma1, Ma2,KS, G3, G4}, due to the simultaneous use of the f\/irst, second and the sophisticated non-linear Spencer complex, and the identif\/ications needed to introduce them.

Section~\ref{partcon} contains  the construction of  partial connections on $J^kV$. These partial connections are fundamental for Section~\ref{formaliso}. Section~\ref{sLieeq} introduces the basic facts on linear Lie equations and the associated groupoids.  Section~\ref{formaliso} presents the def\/inition of formal isomorphism of linear Lie equations, and the proof of Theorem~\ref {iso}. Section~\ref{ila} introduces the def\/inition of intransitive Lie algebras, and the notion of isomorphism of these algebras. The last section discuss the classif\/ication of  intransitive linear Lie equations of order one in the plane, with symbol $g^1$ of dimension one. This classif\/ication contains the examples introduced above.

\section{Preliminaries}\label{pre}
In this section, we present some background material. The main references for this section are~\cite{{Ma1},{Ma2},{KS}}, and we will try
 to maintain the exposition as self-contained as possible, principally introducing geometrical proofs for actions of invertible sections of ${Q}^{k+1}M$ on sections of~$T\oplus J^kTM$.

\subsection{Groupoids and algebroids of jets}

Let be $M$ a manifold and $Q^kM=Q^k$ the manifold of $k$-jets of local dif\/feomorphisms of $M$. This manifold has a natural structure of Lie groupoid given by composition of jets
\[j^k_{f(x)}g.j^k_xf=j^k_x(gf),\]
and inversion
\[\big(j^k_xf\big)^{-1}=j^k_{f(x)}f^{-1},\]
where  $f:U\rightarrow V$, $g:V\rightarrow W$ are local dif\/feomorphisms of $M$, and $x\in U$. The groupoid $Q^k$ has a natural submanifold of identities $I=j^k\mbox{id}$, where $\mbox{id}$ is the identity function of $M$. Then we have a natural identif\/ication of $M$ with $I$, given by $x\mapsto I(x)$. Therefore we can think of $M$ as a submanifold of $Q^k$. There are two submersions  $\alpha,\beta:Q^k\rightarrow M$, the canonical projections \emph{source}, $\alpha(j^k_xf)=x$, and \emph{target}, $\beta(j^k_xf)=f(x)$.  We also consider $\alpha$ and $\beta$ with values in $I$, by the above identif\/ication of $M$ with $I$.

There are natural projections $\pi^k_l=\pi_l:Q^k\rightarrow Q^l$, for $k\geq l\geq 0$, def\/ined by $\pi_l(j^k_xf)=j^l_xf$. Observe that $Q^0=M\times M$ and $\pi^k_0=(\alpha,\beta)$. The projections $\pi^k_l$ commute with the operations of composition and inversion in $Q^k$.

We denote by $Q^k(x)$ the $\alpha$-f\/iber of $Q^k$ on $x\in M$, or $Q^k(x)=\alpha^{-1}(x)$; by $Q^k(\cdot,y)$ the $\beta$-f\/iber of $Q^k$ on $y\in M$, in another way, $Q^k(\cdot,y)=\beta^{-1}(y)$ and $Q^k(x,y)=Q^k(x)\cap Q^k(\cdot,y)$. The set $Q^k(x,x)$ is a group, the so-called \emph{isotropy group} of $Q^k$ at point $x$. If $U, V$ are  open sets of $M$, $Q^k(U)=\cup_{x\in U}Q^k(x)$, $Q^k(\cdot,V)=\cup_{y\in V}Q^k(\cdot,y)$, and $Q^k(U,V)=Q^k(U)\cap Q^k(\cdot,V)$.

A \emph{$($differentiable$)$ section} $F$ of $Q^k$ def\/ined on an open set $U$ of $M$ is a dif\/ferentiable map $F:U\rightarrow Q^k$ such that $\alpha (F(x))=x$. If $\beta (F(U))=V$ and $f=\beta F:U\rightarrow V$ is a dif\/feomorphism, we say that the section~$F$ is \emph{invertible}. We write $U=\alpha(F)$ and $V=\beta(F)$. An invertible section~$F$ of $Q^k$ is said \emph{holonomic} if there exists a dif\/feomorphism $f:\alpha(F)\rightarrow \beta(F)$ such that $F=j^kf$. In this case, $\beta F=f$.

We denote by $\mathcal{Q}^k$ the set of invertible sections of $Q^k$. Naturally $\mathcal{Q}^k$ has a structure of groupoid.  If $F,H\in\mathcal{Q}^k$ with $\beta (F)=\alpha (H)$, then $HF(x)=H(f(x))F(x)$ and $F^{-1}(y)=F(f^{-1}(y))^{-1}$, $y\in \beta (F)$.

Similarly, we can introduce the groupoid of $l$-jets of invertible sections of $Q^k$, and we denote this groupoid by $Q^lQ^k$.
 We have the inclusions
\begin{gather}
\lambda^l: \ Q^{k+l}\rightarrow Q^lQ^k,\nonumber\\
\phantom{\lambda^l:}{} \ \  j^{k+l}_xF\mapsto j^l_xj^kF.\label{lambdal}
\end{gather}
An invertible section $F$, with $\alpha(F)=U$, $\beta(F)=V$, def\/ines a dif\/feomorphism
	\begin{gather*}
	\tilde F: \ Q^k(\cdot,U)\rightarrow  Q^k(\cdot,V),\\
\phantom{\tilde F:}{} \ \ X\mapsto F(\beta(X))X.
	\end{gather*}
The dif\/ferential $\tilde F_*: TQ^k(\cdot,U)\rightarrow TQ^k(\cdot,V)$ depends, for each $X\in Q^k(\cdot,U)$, only on $j^1_{\beta(X)}F$. This def\/ines an action
\begin{gather}
	j^1_{\beta(X)}F: \ T_XQ^k\rightarrow  T_{F(\beta(X))X}Q^k,\nonumber\\
\phantom{j^1_{\beta(X)}F:}{} \ \ v\mapsto j^1_{\beta(X)}F\cdot v=(\tilde F_*)_X(v).\label{jetf}
	\end{gather}
Then~\eqref{jetf} def\/ines a left action of $Q^1Q^k$ on $TQ^k$
\begin{gather}
 Q^1Q^k\times TQ^k\rightarrow  TQ^k,	\nonumber\\
	\big(j^1_{\beta(X)}F,v\in T_XQ^k\big)\mapsto j^1_{\beta(X)}F\cdot v\in T_{F(\beta(X))X}Q^k.\label{actionl}
\end{gather}
If $V_\beta^k\equiv\ker\beta_*\subset TQ^k$ denotes the subvector bundle of $\beta_*$ vertical vectors, then the action~\eqref{jetf} depends only on
$F(\beta (X))$:
\begin{gather}
 Q^k\times V_{\beta}^k\rightarrow  V_{\beta}^k,\nonumber\\
	(F({\beta(X)}),v\in (V_{\beta}^k)_X)\mapsto F({\beta(X)})\cdot v\in (V_{\beta}^k)_{F(\beta(X))X}.\label{restactionl}
\end{gather}
In a similar way, $F$ def\/ines a right action which is a dif\/feomorphism
\begin{gather}
	\bar F: \ Q^k(V)\rightarrow  Q^k(U),\nonumber\\
\phantom{\bar F:}{} \ \	X\mapsto XF(f^{-1}(\alpha(X))).\label{barF}	
\end{gather}
The dif\/ferential $\bar F_*$ of $\bar F$ induces the right action
\begin{gather}
 TQ^k\times Q^1Q^k\rightarrow  TQ^k,\nonumber\\
	\big(v\in T_XQ^k,j^1_{f^{-1}(\alpha(X))}F\big)\mapsto v\cdot j^1_{f^{-1}(\alpha(X))}F=(\bar F_*)_X(v).\label{actionr}
\end{gather}
As $\beta(YX)=\beta(Y)$, it follows that
\begin{equation*}
\beta_*\big(v\cdot j^1_{f^{-1}(\alpha(X))}F\big)=\beta_*(v),
\end{equation*}
where $\beta_*:TQ^k\rightarrow TM$ is the dif\/ferential of $\beta:Q^k\rightarrow M$.
We verify from \eqref{barF} that the function~$\bar F$ restricted to the $\alpha$-f\/iber $Q^k(y)$ depends only on the value of $F$ in $f^{-1}(y)$. If $V_{\alpha}^k=\ker \alpha_*$, then the right action~\eqref{actionr} depends only on the value of $F$ at each point, and the action~(\ref{actionr}) by restriction gives the action
\begin{gather}
V_{\alpha}^k\times Q^k\rightarrow  V_{\alpha}^k,\nonumber\\
	\big(v\in (V_{\alpha}^k)_X,Y\in Q^k(\cdot,\alpha(X))\big)\mapsto v\cdot Y\in (V_{\alpha}^k)_{XY}.\label{actionvert}
	\end{gather}

A vector f\/ield $\bar\xi$
on $Q^k$ with values in $V_{\alpha}^k$ is said \emph{right invariant} if $\bar\xi(XY)=\bar\xi(X)\cdot Y$. The vector f\/ield $\bar\xi$
is determined by its restriction $\xi$ to $I$.

 Let be $T=TM$ the tangent bundle of $M$, and  $\mathcal T$ the sheaf of germs of local sections of $T$. We denote by $J^kT$ the vector bundle of $k$-jets of local sections of $\pi:T\rightarrow M$. Then $J^kT$ is a vector bundle on $M$, and we also denote by $\pi:J^kT\rightarrow M$ the map $\pi(j^k_x\theta)=x$. If $\theta:U\subset M\rightarrow TM$ is a local section, and $f_t$ is the 1-parameter group of local dif\/feomorphisms of $M$ such  that $\frac{d}{dt}f_t|_{t=0}=\theta$, then we get, for $x\in U$,
 \[
 \frac{d}{dt}j^k_xf_t|_{t=0}=j^k_x\theta.
 \]
  This means we have a natural identif\/ication
  \[V_{\alpha}^k|_I=J^kT.\]
   Therefore, as $TQ^k|_I=TI\oplus V_{\alpha}^k|_I$,
   \begin{equation}\label{tqki}TQ^k|_I\cong T\oplus J^kT,
   \end{equation}
 and if we denote by
 \[
 \check{J}^kT=T\oplus J^kT,
 \]
 we have $TQ^k|_I\cong \check{J}^kT$. Observe that $\check{J}^kT $ is a vector bundle on $M$. The restriction of $\beta_*:TQ^k\rightarrow T$ to $TQ^k|_I$, and the isomorphism $TQ^k|_I\cong \check{J}^kT$ def\/ines the map
 \begin{gather}
 \beta_*: \ \check{J}^kT\rightarrow T,\nonumber\\
\phantom{\beta_*:}{} \ \  v+j_x^k\theta\in (\check{J}^kT)_x\mapsto v+\theta(x),\label{beta}
 \end{gather}
which we denote again by $\beta_*$. For more details on this identif\/ication and the map $\beta_*$, see the Appendix of~\cite{KS}, in particular pages~260 and~274.

If $\xi$ is a section of $J^kT$ on $U\subset M$, let be
\[
\bar\xi(X)=\xi(\beta(X))\cdot X
\]
the right invariant vector f\/ield on $Q^k(\cdot,U)$. Then $\bar\xi$ has $\bar F_t$, $-\epsilon<t<\epsilon$, as the 1-parameter group of dif\/feomorphisms induced by invertible sections $F_t$ of $Q^k$ such that
\[
\frac{d}{dt}\bar F_t|_{t=0}=\bar\xi.
\]
Therefore, $F_0=I$ and
\[\frac{d}{dt} F_t(x)|_{t=0}=\xi(x).
\]

\begin{definition}
The vector bundle $J^kT=V_{\alpha}^k|_I$ on $M$ is the \emph{$($differentiable$)$ algebroid}
associated with the groupoid $Q^k$.
\end{definition}

The Lie bracket $\zerobra{\,\,}{\,}_k$ on local sections of $J^kT$  is well def\/ined,  given by
\begin{equation}\label{Liebra}
\zerobra{\xi}{\eta}_k=\left[\bar\xi,\bar\eta\right]|_I,
\end{equation}
where $\xi$, $\eta$ are sections of $\pi^k:J^kT\rightarrow M$ def\/ined on an open set $U$ of $M$ (or $I$).
\begin{proposition}
If $f$ is a real function on $U$, and $\xi$, $\eta$ sections of $J^kT$ on $U$,
then
\[
\zerobra{f\xi}{\eta}_k=f\zerobra{\xi}{\eta}_k-(\beta_*\eta)(f)\xi.
\]
\end{proposition}

\begin{proof} As $\overline {f\xi}=(f\beta)\bar\xi$, it follows
\begin{gather*}
\zerobra{f\xi}{\eta}_k=\left[(f\beta)\bar\xi,\bar\eta\right]|_I
=\left((f\beta)\left[\bar\xi,\bar\eta\right]-\bar\eta(f\beta)\bar\xi\right)|_I
=f\zerobra{\xi}{\eta}_k-(\beta_*\eta)(f)\xi.\tag*{\qed}
\end{gather*}
\renewcommand{\qed}{}
\end{proof}

If $J^k\mathcal{T}$ denotes the sheaf of germs  of local sections of $J^kT$, then $J^k\mathcal{T}$ is a Lie algebra sheaf, with the Lie bracket $\zerobra{\,\,}{\,\,}_k$.

\begin{proposition}\label{bracketk}
The bracket $\zerobra{\,\,}{\,\,}_k$ on ${J}^k\mathcal{T}$ is determined by:
\begin{enumerate}\itemsep=0pt
	\item[$(i)$] $\zerobra{j^k\xi}{j^k\eta}_k=j^k[\xi,\eta]$, $\xi,\eta\in \mathcal{T}$,
	\item[$(ii)$] $\zerobra{\xi_k}{f\eta_k}_k=f\zerobra{\xi_k}{\eta_k}_k+(\beta_*\xi_k)(f)\eta_k$,
\end{enumerate}
where $\xi_k,\eta_k\in J^k\mathcal{T}$ and $f$ is a real function on $M$.
\end{proposition}

 We denote by the same symbols as above the projections $\pi^k_l=\pi_l:J^kT\rightarrow J^lT$, $l\geq 0$, def\/ined by $\pi_l(j^k_x\theta)=j^l_x\theta$. If $\xi_k$ is a point or a section of $J^kT$, let be $\xi_l=\pi^k_l(\xi_k)$. The vector bundle $J^0T$ is isomorphic to $T$ by  the map $\beta_*:J^0T\rightarrow T$, where  $\beta_*(j^0_x\theta)=\theta(x)$, see \eqref{beta}. However, $\beta_*:J^kT\rightarrow T$  is not equal to $\pi^k_0:J^kT\rightarrow J^0T$, but they are isomorphic maps.

Again, we have the canonical inclusions, and we use the same notation as \eqref{lambdal},
\begin{gather*}
\lambda^l: \ J^{k+l}T\rightarrow J^lJ^kT,\nonumber\\
\phantom{\lambda^l:}{} \ \ j^{k+l}_x\theta\mapsto j^l_xj^k\theta\label{lkl}.
\end{gather*}
for $\theta\in\mathcal T$.

Analogously to the def\/inition of holonomic sections of $Q^k$, a section $\xi_k$ of $J^kT$ is \emph{holonomic} if there exists $\xi\in\mathcal{T}$ such that  $\xi_k=j^k\xi$. Therefore, if $\xi_k$ is holonomic, we have $\xi_k=j^k(\beta_*\xi_k)$.

If $\theta:U\subset M\rightarrow J^kT$ is a section, let be $\xi=j^1_x\theta\in J^1_xJ^kT$, $x\in U$. Then $\xi$ can be identif\/ied to the linear application
\begin{gather*}
	\xi: \ T_x\rightarrow T_{\theta(x)}J^kT,\\
\phantom{\xi:}{} \ \	v\mapsto \theta_*(v)	.
\end{gather*}
If $\eta\in J^1_xJ^kT$ is given by $\eta=j^1_x\mu$, with $\mu(x)=\theta(x)$,   then $(\pi)_*(\eta-\xi)v=0$, and we remember that $\pi:J^kT\rightarrow M$ is def\/ined by $\pi(j^k_x\theta )=x$.
So $\eta-\xi\in T_x^*\otimes V_{\pi^1_0(\xi)}J^kT$, where $VJ^kT=\ker\pi_{*}$.   However $J^kT$ is a vector bundle, then $V_{\pi^1_0\xi}J^kT\cong J^k_xT$, so $\eta-\xi\in T_x^*\otimes J^k_xT$. The sequence
\begin{gather}\label{kernelpi}
0\rightarrow T^*\otimes J^kT\rightarrow J^1J^kT\stackrel{\pi^1_0}{\rightarrow}J^kT\rightarrow 0
\end{gather}
obtained in this way is exact, and we get an af\/f\/ine structure on $J^1J^kT$.

The linear operator  $D$  def\/ined by
\begin{gather}
	D: \ J^k\mathcal{T}\rightarrow \mathcal{T}^*\otimes J^{k-1}\mathcal{T},\nonumber\\
\phantom{D:}{} \ \	\xi_k\mapsto D\xi_k=j^1\xi_{k-1}-\lambda^1(\xi_k),\label{D}	
\end{gather}
is the \emph{linear Spencer operator}. We remember that $\xi_{k-1}=\pi^k_{k-1}\xi_k$ and
\begin{gather*}
	\lambda^1: \ J^kT\rightarrow J^1J^{k-1} T,\\
\phantom{\lambda^1:}{} \ \ 	j^k_x\xi\mapsto j^1_x(j^{k-1}\xi)	.
\end{gather*}
The dif\/ference in \eqref{D} is done in $J^1J^{k-1}T$ and is in $T^*\otimes J^{k-1}T$, by \eqref{kernelpi}.

The operator $D$ is null on a section $\xi_k$ if and only if it is holonomic, i.e., $D\xi_k=0$ if and only if there exists $\theta\in\mathcal{T}$ such that $\xi_k=j^k\theta$.
\begin{proposition}\label{propriedadesD}The operator $D$ is characterized by
\begin{enumerate}\itemsep=0pt
	\item[$(i)$] $Dj^k=0$,
	\item[$(ii)$] $D(f\xi_k)=df\otimes \xi_{k-1}+fD\xi_k$,
\end{enumerate}
with $\xi_k\in J^k\mathcal T$, $\xi_{k-1}=\pi_{k-1}\xi_k$ and $f$ is a real function on $M$.
\end{proposition}
For a proof, see \cite{KS}.

The operator $D$ extends to
\begin{gather*}
	D: \ \wedge^l\mathcal{T}^*\otimes J^k\mathcal{T}\rightarrow \wedge^{l+1}\mathcal{T}^*\otimes J^{k-1}\mathcal{T},\nonumber\\
\phantom{D:}{} \ \ \omega\otimes\xi_k\mapsto D(\omega\otimes\xi_k)=d\omega\otimes \xi_{k-1}+(-1)^l\omega\wedge D\xi_k.\label{extD}
	\end{gather*}

\subsection{The calculus on the diagonal}\label{calcdiag}

Next,  following \cite{Ma1,Ma2,KS}, we will relate $\check{J}^k\mathcal{T}$ to vector f\/ields along the diagonal of $M\times M$ and actions of sections in $\mathcal{Q}^k$ to dif\/feomorphisms of $M\times M$ which leave the diagonal invariant.

We denote the \emph{diagonal} of $M\times M$ by $\Delta=\{(x,x)\in M\times M|x\in M\}$, and by $\rho_1:M\times M\rightarrow M$ and $\rho_2:M\times M\rightarrow M$, the f\/irst and second projections, respectively. The restrictions~$\rho_1|_\Delta$ and~$\rho_2|_\Delta$ are dif\/feomorphisms of $\Delta$ on $M$. A sheaf on $M$ will be identif\/ied to its inverse image by~$\rho_1|_\Delta$. For example, if $\mathcal{O}_M$ denotes the sheaf of germs  of real functions on $M$, then we will write~$\mathcal{O}_M$ on~$\Delta$ instead of $(\rho_1|_{\Delta})^{-1}\mathcal{O}_M$. Therefore, a $f\in \mathcal{O}_M$ will be considered in~$\mathcal{O}_{\Delta}$ or in~$\mathcal{O}_{M\times M}$ through the map $f\mapsto f\circ \rho_1$.

We denote by $\mathcal{T}(M\times M)$ the sheaf of germs of local sections of $T(M\times M)\rightarrow M\times M$; by~$\mathcal{R}$ the subsheaf in Lie algebras of $\mathcal{T}({M}\times {M})$, whose elements are vector f\/ields $\rho_1$-projectables; by~$\mathcal{H_R}$ the subsheaf in Lie algebras of $\mathcal{R}$ that projects on $0$ by~$\rho_2$, i.e.\ $\mathcal{H_R}=(\rho_2)_*^{-1}(0)\cap\mathcal{R}$; and by~$\mathcal{V_R}$ the subsheaf in Lie algebras def\/ined by $\mathcal{V_R}=(\rho_1)_*^{-1}(0)\cap\mathcal{R}$. Clearly, \[\mathcal{R}=\mathcal{H_R}\oplus \mathcal{V_R},\] and \[[\mathcal{H_R},\mathcal{V_R}]\subset\mathcal{V_R}.\]
 Then
\[(\rho_1)_*: \ \mathcal{H_R}\widetilde{\longrightarrow}
\mathcal {T}\]
is an isomorphism, so we identify $\mathcal{H_R}$ naturally with $\mathcal{T}$ by this isomorphism, and utilize both notations indistinctly.
\begin{proposition}\label{bili}
The Lie bracket in $\mathcal{R}$ satisfies:
\begin{equation*}
[v+\xi,f(w+\eta)]=v(f)(w+\eta)+f[v+\xi,w+\eta],
\end{equation*}
with $f\in\mathcal{O}_M$, $v,w\in \mathcal{H_R}$, $\xi,\eta\in\mathcal{V_R}$.
In particular, the Lie bracket in $\mathcal{V_R}$ is $\mathcal{O}_M$-bilinear.
\end{proposition}
\begin{proof} Let be $f\in \mathcal{O}_M$, $\xi,\eta\in\mathcal{V_R}$. Then
\[[v+\xi,(f\circ\rho_1)(w+\eta)]=(v+\xi)(f\circ\rho_1)(w+\eta)+(f\circ\rho_1)[v+\xi,w+\eta].\]
As $f\circ\rho_1$ is constant on the submanifolds $\{x\}\times M$ and
$\xi$ is tangent to them, we obtain $\xi(f\circ\rho_1)=0$, and the proposition is proved.
\end{proof}

A vector f\/ield in $\mathcal{V_R}$ is given by a family of sections of $\mathcal{T}$ parameterized by an open set of~$M$.
Therefore there exists a surjective morphism
\begin{gather*}
	\Upsilon_k: \ \mathcal{R}\rightarrow \mathcal{T}\oplus J^k\mathcal{T},\nonumber\\
\phantom{\Upsilon_k:}{} \ \ v+\xi\mapsto v+\xi_k,\label{R}	
\end{gather*}
where $v\in\mathcal{H_R}$, $\xi\in\mathcal{V_R}$, and \[\xi_k(x)=j^k_{(x,x)}(\xi|_{\{x\}\times M}).\]
The kernel of  morphism $\Upsilon_k$ is the subsheaf $\mathcal{V_R}^{k+1}$ of $\mathcal{V_R}$ constituted by vector f\/ields that are null on $\Delta$ at order $k$.  Therefore $\mathcal{R}/\mathcal{V_R}^{k+1}$ is null outside  $\Delta$. It will be considered as a sheaf on $\Delta$, and the sections in the quotient as sections on open sets of $M$. So the sheaf $\mathcal{R}/\mathcal{V_R}^{k+1}$ is isomorphic to the sheaf os germs of sections of the vector bundle $T\oplus J^kT$ on $M$. So
we have the isomorphism of sheaves on $M$,
\[
\mathcal{R}/\mathcal{V_R}^{k+1}\equiv\mathcal{T}\oplus J^k\mathcal{T}.
\]
 We usually denote by $\check J^kT=T\oplus J^kT$ and  $\check J^k\mathcal T=\mathcal T\oplus J^k\mathcal T$.

As \[[\mathcal{R},\mathcal{V_R}^{k+1}]\subset \mathcal{V_R}^k,\]
the bracket
on $\mathcal{R}$ induces a bilinear antisymmetric map, which we call the  \emph{first bracket} of order~$k$,
\begin{equation}\label{jka}
\firstbra{\,\,}{\,\,}_k=(\mathcal{T}\oplus J^k\mathcal{T})\times (\mathcal{T}\oplus J^k\mathcal{T})\rightarrow (\mathcal{T}\oplus J^{k-1}\mathcal{T})
\end{equation}
def\/ined by
\[\firstbra{v+\xi_k}{w+\eta_k}_k=\Upsilon_{k-1}([v+\xi,w+\eta]),
\]
where $\Upsilon_k(\xi)=\xi_k$ and $\Upsilon_k(\eta)=\eta_k$.

It follows from Proposition~\ref{bili} that $\firstbra{\,\,}{\,\,}_k$ satisf\/ies:
\begin{gather}
\firstbra{v+\xi_k}{f(w+\eta_k)}_k=v(f)(w+\eta_{k-1})+f\firstbra{v+\xi_k}{w+\eta_k}_k,\label{noombi}\\
\firstbra{\firstbra{v+\xi_k}{w+\eta_k}_k}{z+\theta_{k-1}}_{k-1}+
\firstbra{\firstbra{w+\eta_k}{z+\theta_k}_k}{v+\xi_{k-1}}_{k-1}\nonumber\\
\phantom{\firstbra{\firstbra{v+\xi_k}{w+\eta_k}_k}{z+\theta_{k-1}}_{k-1}}{}
+\firstbra{\firstbra{z+\theta_k}{v+\xi_k}_k}{w+\eta_{k-1}}_{k-1}=0,
\nonumber
\end{gather}
for $v,w,z\in\mathcal{T}$, $\xi_k,\eta_k,\theta_k\in J^k\mathcal{T}$, $f\in\mathcal{O}_M$. In particular, the f\/irst bracket
is $\mathcal{O}_M$-bilinear on~$J^k\mathcal{T}$. Also,
\[\firstbra{J^0\mathcal{T}}{J^0\mathcal{T}}_0=0.\]

The following proposition relates $\firstbra{\,\,}{\,\,}_k$ to the bracket in $\mathcal{T}$ and the linear Spencer opera\-tor~$D$ in $J^k\mathcal{T}$.
\begin{proposition}\label{colchetek}
Let be $v,w,\theta,\mu\in\mathcal{T}$, $\xi_k,\eta_k\in J^k\mathcal{T}$ and $f\in\mathcal{O}_M$. Then:
\begin{enumerate}\itemsep=0pt
\item [$(i)$] $\firstbra{v}{w}_k=[v,w]$, where the  bracket at right is the bracket in $\mathcal{T}$;
\item[$(ii)$] $\firstbra{v}{\xi_k}_k=i(v)D\xi_k$;
\item[$(iii)$] $\firstbra{j^k\theta}{j^k\mu}_k=j^{k-1}[\theta,\mu]$, where the bracket at right is the bracket in $\mathcal{T}$.
\end{enumerate}
\end{proposition}
\begin{proof} $(i)$ This follows from the identif\/ication of $\mathcal{T}$ with $\mathcal{H_R}$.

$(ii)$ First of all, if $\theta\in\mathcal{T}$, let be $\Theta\in\mathcal{V_R}$ def\/ined by $\Theta(x,y)=\theta(y)$. Then $\Upsilon_k(\Theta)=j^k\theta$. If $v\in\mathcal{H_R}$, then $v$ and $\Theta$ are both $\rho_1$ and $\rho_2$ projectables,  $(\rho_1)_*(\Theta)=0$ and $(\rho_2)_*(v)=0$, so we get $[v,\Theta]=0$. Consequently
\begin{equation}\label{bvt}
\firstbra{v}{j^k\theta}_k=\Upsilon_{k-1}([v,\Theta])=0.
\end{equation}
Also by \eqref{noombi}, we have
\begin{equation}\label{vfx}
\firstbra{v}{f\xi_k}_k=v(f)\xi_{k-1}+f\firstbra{v}{\xi_k}_k.
\end{equation}
As \eqref{bvt} and \eqref{vfx} determine $D$ (cf.\ Proposition~\ref{propriedadesD}), we get $(ii)$.

$(iii)$ Given $\theta,\mu\in\mathcal{T}$, we def\/ine $\Theta,H\in\mathcal{V_R}$ as in $(ii)$, $\Theta(x,y)=\theta(y)$ and $H(x,y)=\mu(y)$. Therefore
\begin{gather*}
\firstbra{j^k\theta}{j^k\mu}_k=\firstbra{\Upsilon_k\Theta}{\Upsilon_kH}_k=\Upsilon_{k-1}([\Theta,H])
=j^{k-1}[\theta,\eta].\tag*{\qed}
\end{gather*}
\renewcommand{\qed}{}
\end{proof}

Let be $\tilde{\mathcal{V}}_\mathcal{R}$ the subsheaf in Lie algebras of $\mathcal{R}$ such that $\tilde\xi\in\tilde{\mathcal{V}}_\mathcal{R}$ if and only if $\tilde\xi$ is tangent to the diagonal $\Delta$.
If $\tilde\xi=\xi_H+\xi\in\tilde{\mathcal{V}}_\mathcal{R}$, $\xi_H\in\mathcal{H}$, $\xi\in\mathcal{V}$, then
 \begin{equation*}
 (\rho_1)_*(\xi_H(x,x))=(\rho_2)_*(\xi(x,x)),
 \end{equation*}
 where $\rho_1,\rho_2:M\times M\rightarrow \Delta$. Consequently, if $\xi_k=\Upsilon_k(\xi)$, then $\xi_H=\beta_*(\xi_k)$, where we remember that  $\beta_*:J^kT\rightarrow T$ is def\/ined in \eqref{beta}. We can also  write $\xi_H=(\rho_2)_*(\xi_k)$, since that $(\rho_2)_*\mathcal{V_R}^{k+1}|_{\Delta}=0$. From now on, $\xi_H$ denotes the horizontal component of $\tilde\xi\in\tilde{\mathcal{V}}_\mathcal{R}$, so $\tilde\xi=\xi_H+\xi$, with $\xi_H\in\mathcal{H}$ and $\xi\in\mathcal{V}$.

 We denote by $J^k\tilde{\mathcal{T}}$ the subsheaf of $\mathcal{T}\oplus J^k\mathcal{T}=\check J^k\mathcal{T}$, whose elements are
 \[\tilde\xi_k=\xi_H+\xi_k,\]
  where $\xi_H=\beta_*(\xi_k)$ or $\xi_H=(\rho_2)_*(\xi_k)$. Therefore $J^k\tilde{\mathcal{T}}$ identif\/ies with $\tilde{\mathcal{V}}_\mathcal{R}/\mathcal{V_R}^{k+1}$, since that $\mathcal{V_R}^{k+1}\subset\tilde{\mathcal{V}}_\mathcal{R}$. As
 \begin{equation}\label{VtV}
 [\tilde{\mathcal{V}}_\mathcal{R},\mathcal{V_R}^{k+1}]\subset \mathcal{V_R}^{k+1},
 \end{equation}
 since the vector f\/ields in $\tilde{\mathcal{V}}_\mathcal{R}$ are tangents to $\Delta$, it follows that the bracket in $\tilde{\mathcal{V}}_\mathcal{R}$ def\/ines a~bilinear antisymmetric map,  called the \emph{second bracket}, by
 \begin{gather}
	\secondbra{\,\,}{\,\,}_k: \ J^k\tilde{\mathcal{T}}\times J^k\tilde{\mathcal{T}}\rightarrow J^k\tilde{\mathcal{T}}, \nonumber\\
\phantom{\secondbra{\,\,}{\,\,}_k:}{} \ \	(\xi_H+\xi_k,\eta_H+\eta_k)\mapsto \Upsilon_k([\tilde\xi,\tilde\eta]),\label{bta} 	
 \end{gather}
 where $\Upsilon_k(\tilde\xi)=\xi_H+\xi_k$ and $\Upsilon_k(\tilde\eta)=\eta_H+\eta_k$.
Unlike the f\/irst bracket~\eqref{jka}, we do not lose one order doing the bracket in $J^k\tilde{\mathcal{T}}$.
The second bracket $\secondbra{\,\,}{\,\,}_k$ is a Lie bracket on $J^k\tilde{\mathcal{T}}$. The  Proposition~\ref{tilde} below relates it to the  bracket $\zerobra{\,\,}{\,\,}_k$, def\/ined in \eqref{Liebra}.

The projection
\begin{gather*}
\nu: \ \mathcal{H_R}\oplus\mathcal{V_R}\rightarrow \mathcal{V_R},\\
\phantom{\nu:}{} \ \ v+\xi\mapsto \xi
\end{gather*}
quotients to
\begin{gather*}
	\nu_k: \ \mathcal{T}\oplus J^k{\mathcal{T}}\rightarrow J^k\mathcal{T}, \\
\phantom{\nu_k:}{} \ \ v+{\xi}_k	\mapsto  \xi_k,\label{nu}
\end{gather*}
and $\nu_k:J^k\tilde{\mathcal{T}}\rightarrow J^k{\mathcal{T}}$ is an isomorphism of vector bundles.
\begin{proposition}\label{tilde}
If $\tilde\xi_k,\tilde\eta_k\in J^k\tilde{\mathcal{T}}$, then
\[\zerobra{\xi_k}{\eta_k}_k=\nu_k(\secondbra{\tilde\xi_k}{\tilde\eta_k}_k),\]
where $\xi_k=\nu_k(\tilde\xi_k)$, $\eta_k=\nu_k(\tilde\eta_k)$.
\end{proposition}
\begin{proof} We will verify properties $(i)$ and $(ii)$ of Proposition \ref{bracketk}. If $\theta,\mu\in\mathcal{T}$, let be $\Theta,H\in\mathcal{V_R}$ as in the proof of Proposition \ref{colchetek}. Then:
\begin{gather*}
(i) \  \nu_k(\secondbra{\theta+j^k\theta}{\mu+j^k\mu}_k)=\nu_k(\Upsilon_k([\theta+\Theta,\mu+H]))
 =\Upsilon_k(\nu([\theta,\mu]+[\Theta,H])) \\
 \phantom{(i) \  \nu_k(\secondbra{\theta+j^k\theta}{\mu+j^k\mu}_k)}{}
=j^k([\theta,\mu])  =\zerobra{j^k\theta}{j^k\mu}_k.\\
(ii) \ \nu_k(\secondbra{\tilde\xi_k}{f\tilde\eta_k}_k)=\nu_k(f\secondbra{\tilde\xi_k}{\tilde\eta_k}_k+\xi_H(f)\tilde\eta_k)
	= f\nu_k(\secondbra{\tilde\xi_k}{\tilde\eta_k}_k)+(\beta_*\xi_k)(f)\eta_k.\tag*{\qed}
\end{gather*}
\renewcommand{\qed}{}
\end{proof}
\begin{corollary}
If $\tilde\xi_k,\tilde\eta_k\in J^k\tilde{\mathcal{T}}$, then
\[\nu_k(\secondbra{\tilde\xi_k}{\tilde\eta_k}_k)=\zerobra{\xi_k}{\eta_k}_k=i(\xi_H)D\eta_{k+1}-i(\eta_H)D\xi_{k+1}+\firstbra{\xi_{k+1}}{\eta_{k+1}}_{k+1},\]
where $\xi_H=\beta_*\xi_k,\,\eta_H=\beta_*\eta_k\in\mathcal{T}$ and $\xi_{k+1},\eta_{k+1}\in J^{k+1}\mathcal{T}$ projects on $\xi_k,\,\eta_k$, respectively.
\end{corollary}
\begin{proof}
It follows from Propositions \ref{colchetek} and \ref{tilde}.
\end{proof}

As a consequence of Proposition \ref{tilde}, we obtain that
\begin{gather}
	\nu_k: \ J^k\tilde{\mathcal{T}}\rightarrow J^k\mathcal{T}, \nonumber\\
\phantom{\nu_k:}{} \ \ \tilde{\xi}_k	\mapsto  \xi_k	\nonumber 
\end{gather}
is an isomorphism of Lie algebras sheaves, where the bracket in $J^k\tilde{\mathcal{T}}$ is the second bracket $\secondbra{\,\,}{\,\,}_k$ as def\/ined in \eqref{bta}, and the bracket in $J^k{\mathcal{T}}$ is the  bracket $\zerobra{\,\,}{\,\,}_k$ as def\/ined in \eqref{Liebra}.

In a similar way, we obtain from \eqref{VtV} that we can def\/ine the \emph{third bracket} as
 \begin{gather*}
	\thirdbra{\,\,}{\,\,}_k: \ J^{k+1}\tilde{\mathcal{T}}\times \check{J}^k{\mathcal{T}}\rightarrow \check J^k{\mathcal{T}}, \nonumber \\
\phantom{\thirdbra{\,\,}{\,\,}_k:}{} \ \ (\xi_H+\xi_{k+1},v+\eta_k)\mapsto \Upsilon_k([\tilde\xi,v+\eta])\label{tb} .	 \end{gather*}
where $\tilde\xi\in\tilde{\mathcal{V}}_{\mathcal R}$, $v+\eta\in\mathcal{R}$,

\begin{proposition}\label{propthird}
The third bracket has the following properties:
\begin{enumerate}\itemsep=0pt
\item[$(i)$] $\thirdbra{f\tilde\xi_{k+1}}{g\check\eta_k}_k=f\xi_H(g)\check\eta_k-v(f)g\tilde\xi_{k}
    +fg\thirdbra{\tilde\xi_{k+1}}{\check\eta_k}_k$;
\item[$(ii)$] $\thirdbra{\tilde\xi_{k}}{\firstbra{\check\eta_k}{\check\theta_k}_k}_{k-1}
    =\firstbra{\thirdbra{\tilde\xi_{k+1}}{\check\eta_k}_k}{\check\theta_k}_k
    +\firstbra{\check\eta_k}{\thirdbra{\tilde\xi_{k+1}}{\check\theta_k}_{k}\,}_k$;
\item[$(iii)$]$\thirdbra{\tilde\xi_{k+1}}{\check\eta_k}_k=\firstbra{\tilde\xi_{k+1}}{\check\eta_{k+1}}_{k+1}$,
\end{enumerate}
where $\tilde\xi_{k+1}=\xi_H+\xi_{k+1}\in \tilde{J}^{k+1}\mathcal T$, $\check\theta_k\in \check J^k\mathcal T$, $\check\eta_{k+1} =v+\eta_{k+1}\in\check J^{k+1}\mathcal T$, $\tilde\xi_{k}=\pi_k(\tilde\xi_{k+1})$, $\check\eta_{k} =\pi_k(\check\eta_{k+1})$.
\end{proposition}
\begin{proof} The proof follows the same lines as the proof of Proposition~\ref{tilde}.
\end{proof}

Let's now verify the relationship  between the action of  dif\/feomorphisms of $M\times M$, which are $\rho_1$-projectable and preserve $\Delta$, on $\mathcal{R}$, and actions \eqref{actionl} and \eqref{actionr} of $Q^1Q^k$
on $TQ^k$. Let be~$\sigma$ a~(local) dif\/feomorphism of $M\times M$ that is $\rho_1$-projectable. Then
\[\sigma(x,y)=(f(x),\Phi(x,y))
,\]
that is, $\sigma$ is def\/ined by $f\in \mbox{Dif\/f } M$, and a function
\begin{gather*}	\phi: \ M\rightarrow \mbox{Dif\/f}\;   M , \\
\phantom{\phi:}{} \ \ x	\mapsto  \phi_x,
\end{gather*}
such that $\phi_x(y)=\Phi(x,y)$. Particularly, when $\sigma(\Delta)=\Delta$, then $ \Phi(x,x)=f(x)$, for all $x\in M$, or \[\phi_x(x)=f(x).\]
As a special case,
\begin{equation*}
(\phi_{f^{-1}(x)})^{-1}(x)=f^{-1}(x).\end{equation*}

Let's denote by $\mathcal{J}$ the set of (local) dif\/feomorphisms of $M\times M$ that are $\rho_1$-projectable and preserve $\Delta$.  We naturally have  the application
\begin{gather}
	\mathcal{J}\rightarrow \mathcal{Q}^k, \nonumber\\
\sigma	\mapsto  \sigma_k\label{sigma},
\end{gather}
where $\sigma_k(x)=j^k_x\phi_x$, $x\in M$. If $\sigma'\in \mathcal{J}$, with $\sigma'=(f',\Phi')$, then
\[	(\sigma'\circ\sigma)(x,y)= \sigma'(f(x),\phi_x(y))
	=  (f'(f(x)),\phi'_{f(x)}(\phi_x(y))
	=((f'\circ f)(x),(\phi'_{f(x)}\circ\phi_x)(y)),
\]
and from this it follows
\[(\sigma'\circ\sigma)_k(x)=j^k_x(\phi'_{f(x)}\circ \phi_x)=j^k_{f(x)}\phi'_{f(x)}.j^k_x\phi_x=\sigma'_k(f(x)).\sigma_k(x)=(\sigma'_k\circ\sigma_k)(x),\]
for each $x\in I$. So \eqref{sigma} is a surjective morphism of groupoids. If $\phi\in\mbox{Dif\/f }M$, let be $\tilde\phi\in \mathcal{J}$ given by
\[\tilde\phi(x,y)=(\phi(x),\phi(y)).\]
It is clear that
\[(\tilde\phi)_k=j^k\phi.\]
It follows from def\/initions of $\mathcal J$ and $\mathcal R$ that the action
\begin{gather}
 \mathcal{J}\times\mathcal{R}\rightarrow \mathcal R,\nonumber\\
(\sigma,v+\xi)\mapsto \sigma_*(v+\xi),\label{JR}
\end{gather}
is well def\/ined.
Then $\mathcal V$ and $\tilde{\mathcal V}$ are invariants by the action of $\mathcal J$.


\begin{proposition}\label{actionsigma}
Let be $\sigma\in\mathcal J$, $v\in\mathcal {H_R}$, $\xi\in\mathcal {V_R}$. We have:
\begin{enumerate}\itemsep=0pt
\item [$(i)$] $(\sigma*\xi)_k=\lambda^1\sigma_{k+1}.\xi_k.\sigma_k^{-1}$;
\item [$(ii)$]$(\sigma*\tilde\xi)_k=f_*(\xi_H)+j^1\sigma_{k}.\xi_k.\sigma_k^{-1}$ ;
\item [$(iii)$] $(\sigma_*v)_k=f_*(v)+(j^1\sigma_k.v.\lambda^1\sigma_{k+1}^{-1}-\lambda^1\sigma_{k+1}.v.\lambda^1\sigma_{k+1}^{-1}).$
\end{enumerate}
\end{proposition}

\begin{proof} If
\[\sigma(x,y)=(f(x),\phi_x(y)),\]
then
\[\sigma^{-1}(x,y)=\big(f^{-1}(x),\big(\phi^{-1}\big)_x(y)\big),\]
where
\[\big(\phi^{-1}\big)_x=\big(\phi_{f^{-1}(x)}\big)^{-1}.\]

$(i)$ Let be
\[\xi=\frac{d}{dt}V_t\big|_{t=0},\]
 where $V_t(x,y)=(x,\eta^t_x(y))$, with $\eta^t_x\in\mbox{Dif\/f }M$ for each $t$, and  $g_t(x)=\eta^t_x(x)$. Then
\[\big(\sigma\circ V_t\circ \sigma^{-1}\big)(x,y)=\big(x,\big(\phi_{f^{-1}(x)}\circ \eta^t_{f^{-1}(x)}\circ (\phi_{f^{-1}(x)})^{-1}\big)(y)\big),\]
and
\[(\sigma_*\xi)(x,y)=\frac{d}{dt}\big(\phi_{f^{-1}(x)}\circ \eta^t_{f^{-1}(x)}\circ (\phi_{f^{-1}(x)})^{-1}\big)\big|_{t=0}(y).
\]
Consequently
\begin{gather*}
\Upsilon_k(\sigma_*\xi)(x)=j^k_x\left(\frac{d}{dt}\big(\phi_{f^{-1}(x)}\circ \eta^t_{f^{-1}(x)}\circ (\phi_{f^{-1}(x)})^{-1}\big)|_{t=0}\right) \\
\phantom{\Upsilon_k(\sigma_*\xi)(x)}{}=\frac{d}{dt}\big(j^k_{(g_t\circ f^{-1})(x)}\phi_{f^{-1}(x)}\circ j^k_{f^{-1}(x)}\eta^t_{f^{-1}(x)}\circ j^k_x(\phi_{f^{-1}(x)})^{-1}\big)\big|_{t=0} \\
\phantom{\Upsilon_k(\sigma_*\xi)(x)}{}=\frac{d}{dt}\big((\tilde\phi_{f^{-1}(x)})_k((g_t\circ f^{-1})(x)). j^k_{f^{-1}(x)}\eta^t_{f^{-1}(x)}.((\tilde\phi_{f^{-1}(x)})^{-1})_k(x)\big)\big|_{t=0} \\
\phantom{\Upsilon_k(\sigma_*\xi)(x)}{}=\frac{d}{dt}\big(j^1_{f^{-1}(x)}(\tilde\phi_{f^{-1}(x)})_k. j^k_{f^{-1}(x)}\eta^t_{f^{-1}(x)}.((\tilde\phi_{f^{-1}(x)})^{-1})_k(x)\big)|_{t=0} \\
\phantom{\Upsilon_k(\sigma_*\xi)(x)}{}=\lambda^1(\sigma_{k+1}(f^{-1}(x))).\xi_k(f^{-1}(x)).\sigma_k^{-1}(x),
\end{gather*}
since that
\[j^1_{f^{-1}(x)}(\tilde\phi_{f^{-1}(x)})_k=j^1_{f^{-1}(x)}j^k\phi_{f^{-1}(x)}=\lambda^1(j^{k+1}_{f^{-1}(x)}\phi_{f^{-1}(x)})=\lambda^1(\sigma_{k+1}(f^{-1}(x))).\]
So we proved
\[(\sigma*\xi)_k=\lambda^1\sigma_{k+1}.\xi_k.\sigma_k^{-1}.\]

$(ii)$ Let be, as in $(i)$,  $\xi=\frac{d}{dt}V_t|_{t=0}$, where $V_t(x,y)=(x,\eta^t_x(y))$, with $\eta^t_x\in\mbox{Dif\/f }M$ for each $t$, and  $g_t(x)=\eta^t_x(x)$. Then
\[\tilde \xi=\frac{d}{dt}\tilde V_t\big|_{t=0},\]
 where $\tilde V_t(x,y)=(g_t(x),\eta^t_x(y))$.
Therefore
\[\big(\sigma\circ\tilde V_t\circ \sigma^{-1}\big)(x,y)=
\big(f\circ g_t\circ f^{-1}(x),\big(\phi_{g_t\circ f^{-1}(x)}\circ \eta^t_{f^{-1}(x)}\circ (\phi_{f^{-1}(x)})^{-1}\big)(y)\big),\]
and
\[
(\sigma_*\tilde\xi)(x,y)=f_*\xi_H(x)+\frac{d}{dt}\big(\phi_{g_t\circ f^{-1}(x)}\circ \eta^t_{f^{-1}(x)}\circ (\phi_{f^{-1}(x)})^{-1}\big)\big|_{t=0}(y).
\]
By projecting, we obtain
\begin{gather*}
\Upsilon_k(\sigma_*\tilde\xi)(x)=f_*\xi_H(x)+j^k_x
\left(\frac{d}{dt}\big(\phi_{g_t\circ f^{-1}(x)}\circ \eta^t_{f^{-1}(x)}\circ (\phi_{f^{-1}(x)})^{-1}\big)|_{t=0}\right) \\
\phantom{\Upsilon_k(\sigma_*\tilde\xi)(x)}{}
=f_*\xi_H(x)+\frac{d}{dt}\big(j^k_{(g_t\circ f^{-1})(x)}\phi_{g_t\circ f^{-1}(x)}\circ j^k_{f^{-1}(x)}\eta^t_{f^{-1}(x)}\circ j^k_x(\phi_{f^{-1}(x)})^{-1}\big)\big|_{t=0} \\
\phantom{\Upsilon_k(\sigma_*\tilde\xi)(x)}{}
=f_*\xi_H(x)+\frac{d}{dt}\big(\sigma_k((g_t\circ f^{-1})(x)). j^k_{f^{-1}(x)}\eta^t_{f^{-1}(x)}.(\sigma_k)^{-1}(x)\big)\big|_{t=0} \\
\phantom{\Upsilon_k(\sigma_*\tilde\xi)(x)}{}
=f_*\xi_H(x)+j^1_{f^{-1}(x)}\sigma_{k}.\xi_k(f^{-1}(x)).\sigma_k^{-1}(x),
\end{gather*}
so
\[(\sigma_*\tilde\xi)_k=f_*\xi_H+j^1\sigma_{k}.\xi_k.\sigma_k^{-1}.\]
Observe that this formula depends only on $\sigma_k$.

$(iii)$ By combining $(i)$ and $(ii)$, we obtain
\[(\sigma*\xi_H)_k=(\sigma_*\tilde\xi)_k-(\sigma*\xi)_k
=(f_*\xi_H+j^1\sigma_{k}.\xi_k.\sigma_k^{-1})-(\lambda^1\sigma_{k+1}.\xi_k.\sigma_k^{-1}).\]
As
\begin{gather*}
j^1\sigma_{k}.\xi_k.\sigma_k^{-1}=j^1\sigma_{k}.\xi_k.\lambda^1\sigma_{k+1}^{-1}
=j^1\sigma_{k}.(\xi_k-\xi_H).\lambda^1\sigma_{k+1}^{-1}+j^1\sigma_{k}.\xi_H.\lambda^1\sigma_{k+1}^{-1}\\
\phantom{j^1\sigma_{k}.\xi_k.\sigma_k^{-1}}{}
=\lambda^1\sigma_{k+1}.(\xi_k-\xi_H).\lambda^1\sigma_{k+1}^{-1}
+j^1\sigma_{k}.\xi_H.\lambda^1\sigma_{k+1}^{-1} \\
\phantom{j^1\sigma_{k}.\xi_k.\sigma_k^{-1}}{}
=\lambda^1\sigma_{k+1}.\xi_k.\sigma_{k}^{-1}-\lambda^1\sigma_{k+1}.\xi_H.\lambda^1\sigma_{k+1}^{-1}
+j^1\sigma_{k}.\xi_H.\lambda^1\sigma_{k+1}^{-1},
\end{gather*}
where  $j^1\sigma_k.(\xi_k-\xi_H)=\lambda^1\sigma_{k+1}.(\xi_k-\xi_H)$ follows from $\beta_*(\xi_k-\xi_H)=0$ by~\eqref{restactionl}.
By replacing this equality above we get
\begin{gather*}
(\sigma*\xi_H)_k =\big(f_*\xi_H+\big(\lambda^1\sigma_{k+1}.\xi_k.\sigma_{k}^{-1}
-\lambda^1\sigma_{k+1}.\xi_H.\lambda^1\sigma_{k+1}^{-1}+j^1\sigma_{k}.\xi_H.\lambda^1\sigma_{k+1}^{-1}\big)\big)\\
\phantom{(\sigma*\xi_H)_k =}{}
-\big(\lambda^1\sigma_{k+1}.\xi_k.\sigma_k^{-1}\big)
=f_*\xi_H-\lambda^1\sigma_{k+1}.\xi_H.\lambda^1\sigma_{k+1}^{-1}+j^1\sigma_{k}.\xi_H.\lambda^1\sigma_{k+1}^{-1}.
\tag*{\qed}
\end{gather*}
  \renewcommand{\qed}{}
\end{proof}

It follows from Proposition~\ref{actionsigma}  that action \eqref{JR} projects on an action $(\,\,)_*$:
\begin{gather}
\mathcal Q^{k+1}\times(\mathcal T\oplus J^k\mathcal T)\rightarrow \mathcal T\oplus J^k\mathcal T,\nonumber\\
 (\sigma_{k+1},v+\xi_k)\mapsto  (\sigma_{k+1})_*(v+\xi_k),\label{GTG}
\end{gather}
where
\[(\sigma_{k+1})_*(v+\xi_k)=f_*v+
\big(j^1\sigma_k.v.\lambda^1\sigma_{k+1}^{-1}-\lambda^1\sigma_{k+1}.v.\lambda^1\sigma_{k+1}^{-1}\big)
+\big(\lambda^1\sigma_{k+1}.\xi_k.\sigma_k^{-1}\big).\]
This action verif\/ies
\begin{equation}\label{colchactions}
\firstbra{(\sigma_{k+1})_*(v+\xi_k)}{(\sigma_{k+1})_*(w+\eta_k)}_k=(\sigma_k)_*(\firstbra{v+\xi_k}{w+\eta_k}_k).
\end{equation}

It follows from Proposition \ref{actionsigma} and~\eqref{GTG} that  $(\sigma_{k+1})_*(\xi_k)(x)$ depends only on the value of~$\sigma_{k+1}(x)$ at the point $x$ where $\xi$ is def\/ined, and $(\sigma_{k+1})_*(v)(x)$ depends on the value of $\sigma_{k+1}$ on a~curve tangent to $v(x)$.

Item $(ii)$ of Proposition~\ref{actionsigma} says that restriction to $J^k\tilde{\mathcal T}$ of action \eqref {GTG} depends only on the section $\sigma_k$, so the action
\begin{gather*}
\mathcal Q^{k}\times J^k\tilde{\mathcal T}\rightarrow  J^k\tilde{\mathcal T},\nonumber\\
 (\sigma_{k},\tilde\xi_k)\mapsto  (\sigma_{k})_*(\tilde\xi_k),
\end{gather*}
is well def\/ined, where
\begin{equation*}
(\sigma_k)_*(\tilde\xi_k)=f_*\xi_H+j^1\sigma_k.\xi_k.\sigma_k^{-1}.
\end{equation*}
In this case, we get
\[\secondbra{(\sigma_k)_*\tilde\xi_k}{(\sigma_k)_*\tilde\eta_k}_k=(\sigma_k)_*\secondbra{\tilde\xi_k}{\tilde\eta_k}_k,
\]
and each $\nu_k(\sigma_k)_*\nu_k^{-1}$ acts as an automorphism of the Lie algebra sheaf $J^k{\mathcal T}$:
\[\zerobra{\nu_k(\sigma_k)_*\nu_k^{-1}\xi_k}{\nu_k(\sigma_k)_*\nu_k^{-1}\eta_k}_k
=\nu_k(\sigma_k)_*\nu_k^{-1}\zerobra{\xi_k}{\eta_k}_k.
\]
If $M$ and $M'$ are two manifolds of the same dimension, we can def\/ine \[Q^k(M,M')=\{j^k_xf:f:U\subset M\rightarrow U'\subset M'\,\mbox{ is a dif\/feomorphism}, \ x\in U\}.\]
The groupoid $Q^k$ acts by the right on $Q^k(M,M')$, and $Q^{'k}=Q^k(M',M')$ acts by the left. Redoing the calculus of this subsection in this context, we obtain the analogous action of \eqref{GTG}:
\begin{gather}
\mathcal Q^{k+1}(M,M')\times(\mathcal T\oplus J^k\mathcal T)\rightarrow \mathcal T'\oplus J^k\mathcal T',\nonumber\\
 (\sigma_{k+1},v+\xi_k)\mapsto  (\sigma_{k+1})_*(v+\xi_k)\label{GTGduplo},
\end{gather}
where $\mathcal Q^{k+1}(M,M')$ denotes the set of invertible sections of $ Q^{k+1}(M,M')$,
\begin{equation*}
(\sigma_{k+1})_*(v+\xi_k)=f_*v+\big(j^1\sigma_k.v.\lambda^1\sigma_{k+1}^{-1}
-\lambda^1\sigma_{k+1}.v.\lambda^1\sigma_{k+1}^{-1}\big)
+\big(\lambda^1\sigma_{k+1}.\xi_k.\sigma_k^{-1}\big),
\end{equation*}
and $f=\beta\sigma^{k+1}:M\rightarrow M'$.  This action also verif\/ies \eqref{colchactions}.

\subsection[The Lie algebra sheaf $\wedge( \check{J}^\infty\mathcal{T})^*\otimes ( \check{J}^\infty\mathcal{T})$]{The Lie algebra sheaf $\boldsymbol{\wedge( \check{J}^\infty\mathcal{T})^*\otimes ( \check{J}^\infty\mathcal{T})}$} \label{las}

In this subsection, we continue to follow the presentation of \cite{Ma1,Ma2}.
We denote by $J^\infty\mathcal{T} $ the projective limit of $J^k\mathcal{T}$, say,
\[J^\infty\mathcal T=\lim\mbox{proj } J^k\mathcal T,\]
and
\[\check{J}^\infty\mathcal T=\mathcal T\oplus J^\infty\mathcal T=\lim\mbox{proj } \check J^k\mathcal T.\]
As $T\oplus J^kT\cong TQ^k|_I$, we have the identif\/ication of $\check{J}^\infty T$ with $\lim\mbox{proj } \Gamma(TQ^k|_I)$, where $\Gamma(TQ^k|_I)$ denotes the sheaf of germs  of local sections of the vector bundle $TQ^k|_I\rightarrow I$.
From the fact that $\mathcal T\oplus J^k\mathcal T$ is a $\mathcal{O}_M$-module, we get  $\check{J}^\infty\mathcal T$ is a $\mathcal{O}_M$-module. In the following, we use the notation
\[\check{\xi}=v+\lim\mbox{proj } \xi_k,\,\, \check{\eta}=w+\lim\mbox{proj } \eta_k\in\mathcal T\oplus J^\infty\mathcal T.\]
We def\/ine the f\/irst bracket in $\check{J}^\infty\mathcal T$ as:
\begin{equation}\label{checkcolch}
\firstbra{\check{\xi}}{\check{\eta}}_\infty=\lim\mbox{proj }\firstbra{v+\xi_k}{w+\eta_k}_k.
\end{equation}
With the bracket def\/ined by \eqref{checkcolch}, $\check J^{\infty} \mathcal T$ is a
\emph{Lie algebra sheaf}. Furthermore,
\begin{equation*}
\firstbra{\check\xi}{f\check \eta}_\infty=v(f)\check \eta+f\firstbra{\check\xi}{\check\eta}_\infty.
\end{equation*}
We extend now, as in \cite{{Ma1},{Ma2},{KS}}, the bracket on $\check{J}^\infty \mathcal T$ to a Nijenhuis bracket (see \cite{FN}) on $\wedge(\check J^\infty \mathcal T)^*\otimes (\check J^\infty \mathcal T)$, where
\begin{equation*}
(\check J^\infty \mathcal T)^*=\lim\mbox{ind  }(\check J^k\mathcal T)^*.
\end{equation*}
We introduce the exterior dif\/ferential $d$ on $\wedge(\check J^\infty \mathcal T)^*$, by:

$(i)$ if $f\in\mathcal O_M$, then $df\in(\check{J}^\infty \mathcal T)^*$ is def\/ined by
\begin{equation*}
\langle df,\check\xi\rangle=v(f).
\end{equation*}

$(ii)$ if $\omega\in(\check{J}^\infty \mathcal T)^*$, then $d\omega\in\wedge^{2}(\check{J}^\infty \mathcal T)^*$ is def\/ined by
\begin{equation*}
\langle d\omega,\check\xi\wedge\check\eta\rangle =\theta(\check\xi)\langle\omega,\check\eta\rangle
-\theta(\check\eta)\langle\omega,\check\xi\rangle
-\langle \omega,\firstbra{\check\xi}{\check\eta}_\infty\rangle,
\end{equation*}
where $\theta(\check\xi)f=\langle df,\check\xi\rangle$.

We extend this operator to forms of any degree as a derivation of degree $+1$
\[d: \ \wedge^r(\check{J}^\infty \mathcal T)^*\rightarrow \wedge^{r+1}(\check{J}^\infty \mathcal T)^*.\]
The exterior dif\/ferential $d$ is linear,
\[d(\omega\wedge\tau)=d\omega\wedge\tau+(-1)^{r}\omega\wedge d\tau,\]
for $\omega\in\wedge^r(\check{J}^\infty \mathcal T)^*$, and $d^2=0$.

Remember that  $(\rho_1)_*:\mathcal T\oplus J^\infty\mathcal T\rightarrow \mathcal T$  is the projection given by the decomposition in direct sum of $\check J^\infty\mathcal T=\mathcal{T}\oplus J^\infty\mathcal T$. (We could use, instead of $(\rho_1)_*$, the natural map $\alpha_*:T\oplus J^kT\rightarrow T$, given by $\alpha_*:TQ^k|_I\rightarrow T$, and the identif\/ication~\eqref{tqki}). Then $(\rho_1)^*:\mathcal T^*\rightarrow (\check{J}^\infty\mathcal T)^* $, and this map extends to $(\rho_1)^*:\wedge\mathcal T^*\rightarrow \wedge(\check{J}^\infty\mathcal T)^* $. If $\omega\in\wedge^r(\check{J}^\infty\mathcal T)^*$, then
\[
\langle (\rho_1)^*\omega,\check\xi_1\wedge\cdots\wedge\check\xi_r\rangle
=\langle \omega,v_1\wedge\cdots \wedge v_r\rangle,\]
where $\check\xi_j=v_j+\xi_j$, $j=1,\dots, r$.
It follows that $d((\rho_1)^*\omega)=(\rho_1)^{*}(d\omega)$.
We identify $\wedge\mathcal T^*$ with its image in
$\wedge(\check{J}^\infty\mathcal T)^*$ by $(\rho_1)^*$, and we write simply $\omega$ instead of $(\rho_1)^*\omega$.

Let be $\textbf{u}=\omega\otimes\check\xi\in\wedge(\check{J}^\infty\mathcal T)^*\otimes (\check{J}^\infty\mathcal T)$, $\tau\in\wedge (\check{J}^\infty\mathcal T)^*$, with $\deg\omega=r$ and $\deg\tau=s$. We also def\/ine $\deg\textbf u =r$. Then we def\/ine the derivation of degree $(r-1)$
\[i(\textbf {u}): \ \wedge^s(\check{J}^\infty\mathcal T)^*\rightarrow \wedge^{s+r-1}(\check{J}^\infty\mathcal T)^*\]
by
\begin{equation}\label{i(u)}
i(\textbf {u})\tau=i(\omega\otimes\check\xi)\tau=\omega\wedge i(\check\xi)\tau
\end{equation}
and the Lie derivative
\[\theta(\textbf {u}): \ \wedge^r(\check{J}^\infty\mathcal T)^*\rightarrow\wedge^{r+s}(\check{J}^\infty\mathcal T)^*\]
by
\begin{equation}\label{theta(u)}
\theta(\textbf {u})\tau=i(\textbf {u})d\tau+(-1)^rd(i(\textbf {u})\tau),
\end{equation}
which is a derivation of degree $r$. If $\textbf {v}=\tau\otimes\check\eta$,
we def\/ine
\begin{equation}\label{buv}
[\textbf {u},\textbf {v}]=[\omega\otimes\check\xi,\tau\otimes\check\eta]
=\omega\wedge\tau\otimes\firstbra{\check\xi}{\check\eta}_\infty
+\theta(\omega\otimes\check\xi)\tau\otimes\check\eta-(-1)^{rs}\theta(\tau\otimes\check\eta)\omega\otimes\check\xi.
\end{equation}
A straightforward calculation shows that:
\begin{equation}\label{buvmod}
[\textbf {u},\tau\otimes\check\eta]=\theta(\textbf{u})\tau\otimes\check\eta+(-1)^{rs}\tau\wedge[\textbf u,\check\eta]-(-1)^{rs+s}d\tau\wedge i(\check\eta)\textbf{u},
\end{equation}
where $i(\check \eta)\textbf u=i(\check\eta)(\omega\otimes\check\xi)=i(\check\eta)\omega\otimes\check\xi$. We verify that
\[[\textbf {u},\textbf v]=-(-1)^{rs}[\textbf {v},\textbf u]
\]
and
\begin{equation}\label{brauv}
[\textbf {u},[\textbf v,\textbf w]]=[[\textbf u,\textbf v],\textbf w]+(-1)^{rs}[\textbf v,[\textbf u,\textbf w]],
\end{equation}
where $\deg\textbf u =r$, $\deg\textbf v =s$.

With this bracket, $\wedge (\check{J}^\infty\mathcal T)^*\otimes (\check{J}^\infty\mathcal T)$ is a \emph{Lie algebra sheaf}.
Furthermore, if
\begin{equation}\label{dtutv}
[\theta(\textbf u),\theta(\textbf v)]=\theta(\textbf u)\theta(\textbf v)-(-1)^{rs}\theta(\textbf v)\theta(\textbf u)
\end{equation}
then
\begin{equation}\label{tutv}
[\theta(\textbf u),\theta(\textbf v)]=\theta([\textbf u,\textbf v]).
\end{equation}

In particular, we have the following formulas:
\begin{proposition}\label{p31}
If $\mbox{\bf u}  , \mbox{\bf v} \in \wedge^1(\check{J}^\infty\mathcal T)^*\otimes (\check{J}^\infty\mathcal T)$, $\omega\in\wedge^1(\check{J}^\infty\mathcal T)^*$, $\check\xi,\check\eta\in \check{J}^\infty\mathcal T$, then:
\begin{gather*}
 (i) \ \ \ \langle \theta(\mbox{\bf u})\omega,\check\xi\wedge\check\eta\rangle =
 \theta(i(\check\xi)\mbox{\bf u})\langle \omega,\check\eta\rangle
 -\theta(i(\check\eta)\mbox{\bf u})\langle \omega,\check\xi\rangle \\
\phantom{(i) \ \ \ \langle \theta(\mbox{\bf u})\omega,\check\xi\wedge\check\eta\rangle =}{}
-\langle \omega,\firstbra{i(\check\xi)\mbox{\bf u}}{\check\eta}_\infty+\firstbra{\check\xi}{i(\check\eta)\mbox{\bf u}}_\infty-i(\firstbra{\check\xi}{\check\eta}_\infty)\mbox{\bf u}\rangle,\\
 (ii) \ \ i(\check\xi)[\mbox{\bf u},\check\eta]=\firstbra{i(\check\xi)\mbox{\bf u}}{\check\eta}_\infty-i(\firstbra{\check\xi}{\check\eta}_\infty)\mbox{\bf u}, \\
(iii) \ \langle [\mbox{\bf u},\mbox{\bf v}],\check\xi\wedge\check\eta\rangle =
\firstbra{i(\check\xi)\mbox{\bf u}}{i(\check\eta)\mbox{\bf v}}_\infty-\firstbra{i(\check\eta)\mbox{\bf u}}{i(\check\xi)\mbox{\bf v}}_\infty \\
\phantom{(iii) \ \langle [\mbox{\bf u},\mbox{\bf v}],\check\xi\wedge\check\eta\rangle =}{}
-i(\firstbra{i(\check\xi)\mbox{\bf u}}{\check\eta}_\infty-\firstbra{i(\check\eta)\mbox{\bf u}}{\check\xi}_\infty-i(\firstbra{\check\xi}{\check\eta}_\infty)\mbox{\bf u})\mbox{\bf v}\\
\phantom{(iii) \ \langle [\mbox{\bf u},\mbox{\bf v}],\check\xi\wedge\check\eta\rangle =}{} -i(\firstbra{i(\check\xi)\mbox{\bf v}}{\check\eta}_\infty-\firstbra{i(\check\eta)\mbox{\bf v}}{\check\xi}_\infty-i(\firstbra{\check\xi}{\check\eta}_\infty)\mbox{\bf v})\mbox{\bf u}.
\end{gather*}
\end{proposition}

\begin{proof} It is a straightforward calculus applying the def\/initions.
\end{proof}

If we def\/ine the groupoid
\begin{gather*}
\mathcal Q^{\infty}=\lim\mbox{proj }\mathcal Q^k,
\end{gather*}
then for $\sigma=\lim\mbox{proj }\sigma_k\in \mathcal Q^\infty$, we obtain, from \eqref{GTG},
\begin{gather*}
\sigma_*: \ \check J^\infty{\mathcal T} \rightarrow  \check J^\infty{\mathcal T}, \\
\phantom{\sigma_*:}{} \ \ \xi=v+\lim\mbox{proj }\xi_k \mapsto \sigma_*\xi=\lim\mbox{proj }(\sigma_{k+1})_*(v+\xi_k),
\end{gather*}
so the action
\begin{gather*}
\mathcal Q^\infty\times\check J^\infty{\mathcal T}\rightarrow\check J^\infty{\mathcal T}, \nonumber\\
(\sigma,\xi)\mapsto\sigma_*\xi,
\end{gather*}
is well def\/ined.
It follows from \eqref{colchactions} that $\sigma_*:\check J^\infty{\mathcal T}\rightarrow \check J^\infty{\mathcal T}$ is an \emph{automorphism of  Lie algebra sheaf}.

Given $\sigma\in\mathcal Q^{\infty}$, $\sigma$ acts on $\wedge (\check J\mathcal T)^*$:
\begin{gather*}
\sigma^*: \ \wedge (\check J\mathcal T)^*\rightarrow\wedge (\check J\mathcal T)^*,\\
\phantom{\sigma^*:}{} \ \ \omega\mapsto\sigma^*\omega,
\end{gather*}
where, if $\omega$ is a $r$-form,
\begin{equation}\label{sigmaomega}
\langle \sigma^*\omega,\check\xi_1\wedge\cdots\wedge\check\xi_r\rangle =\langle \omega, \sigma_*^{-1}(\check\xi_1)\wedge\cdots\wedge\sigma_*^{-1}(\check\xi_r)\rangle.
\end{equation}
Consequently, $\mathcal Q^{\infty}$ acts on $\wedge(\check J^{\infty}\mathcal T)^*\otimes (\check J\mathcal T)$:
\begin{gather*}
\mathcal Q^{\infty}\times(\wedge(\check J^{\infty}\mathcal T)^*\otimes (\check J^{\infty}\mathcal T))\rightarrow\wedge(\check J^{\infty}\mathcal T)^*\otimes (\check J^{\infty}\mathcal T),\nonumber\\
 (\sigma,\textbf u)\mapsto\sigma_*\textbf u,
\end{gather*}
where
\begin{equation}\label{sigmastar}
\sigma_*\textbf u=\sigma_*(\omega\otimes\check\xi)=\sigma^*(\omega)\otimes\sigma_*(\check\xi).
\end{equation}
The action of $\sigma_*$ is an \emph{automorphism  of the Lie algebra sheaf} $\wedge(\check J^{\infty}\mathcal T)^*\otimes (\check J^{\infty}\mathcal T)$, i.e.,
\[[\sigma_*\textbf u,\sigma_*\textbf v]=\sigma_*[\textbf u,\textbf v].\]

\subsection{The f\/irst non-linear Spencer complex}\label{fnsc}
In this subsection we will study the subsheaf $\wedge\mathcal T^*\otimes J^{\infty}\mathcal T$ and introduce linear and non-linear Spencer complexes. Principal references are \cite{{Ma1},{Ma2},{KS}}.
\begin{proposition}\label{p41}
The sheaf $\wedge\mathcal T^*\otimes J^{\infty}\mathcal T$ is a Lie algebra  subsheaf of $\wedge(\check J^{\infty}\mathcal T)^*\otimes (\check J^{\infty}\mathcal T)$, and
\[[\omega\otimes\xi,\tau\otimes \eta]=\omega\wedge\tau\otimes\firstbra{\xi}{\eta}_\infty,\]
where $\omega, \tau\in\wedge\mathcal T^*$, $\xi, \eta\in J^{\infty}\mathcal T$.
\end{proposition}
\begin{proof} Let be $\textbf u=\omega\otimes\xi\in\wedge\mathcal T^*\otimes J^{\infty}\mathcal T$. For any $\tau\in\wedge\mathcal T^*$,  $i(\xi)\tau=0$, then, by applying \eqref{i(u)} we obtain $i(\textbf u)\tau=0$, and by \eqref{theta(u)}, $\theta(\textbf u)\tau=0$. So \eqref{buv}
implies $[\omega\otimes\xi,\tau\otimes \eta]=\omega\wedge\tau\otimes\firstbra{\xi}{\eta}_\infty$.
\end{proof}

Let be the \emph{fundamental form}
\[\chi\in (\check J^\infty\mathcal T)^*\otimes (\check J^\infty\mathcal T)\]
def\/ined by
\[i(\check\xi)\chi=(\rho_1)_*(\check\xi)=v,\]
where $\check\xi=v+\xi\in \mathcal T\oplus J^\infty \mathcal T$. In another words, $\chi$ is the projection of $\check J^\infty\mathcal T$ on $\mathcal T$, parallel to~$J^\infty\mathcal T$.

If $\textbf u=\lim \textbf u_k$, we def\/ine $D\textbf u=\lim D\textbf u_k$.
\begin{proposition}\label{p42}
If $\omega\in\wedge\mathcal T^*$, and $\mbox{\bf u}\in\wedge\mathcal T^*\otimes J^\infty\mathcal T$, then:
\begin{enumerate}\itemsep=0pt
\item [$(i)$] $\theta(\chi)\omega=d\omega$;
\item [$(ii)$] $[\chi,\chi]=0$;
\item [$(iii)$] $[\chi,\mbox{\bf u}]=D\mbox{\bf u}$.
\end{enumerate}
\end{proposition}
\begin{proof} Let be $\check\xi=v+\xi,\, \check\eta=w+\eta\in\mathcal T\oplus J^\infty \mathcal T$.

$(i)$ As $\theta(\chi)$ is a derivation of degree~1, it is enough to prove $(i)$ for $0$-forms $f$ and $1$-forms $\omega\in(\check J^\infty\mathcal T)^*$. From \eqref{theta(u)} we have $\theta(\chi)f=i(\chi)df=df$. It follows from Proposition~\ref{p31}$(i)$ that
\begin{gather*}
\langle \theta(\chi)\omega,\check\xi\wedge\check\eta\rangle= \theta(v)\langle \omega,\check\eta\rangle- \theta(w)\langle\omega,\check\xi\rangle- \langle\omega,\firstbra{v}{\check\eta}_\infty+\firstbra{\check\xi}{w}_\infty
-\chi(\firstbra{\check\xi}{\check\eta}_\infty)\rangle\\
\phantom{\langle \theta(\chi)\omega,\check\xi\wedge\check\eta\rangle}{}
=\theta(v)\langle\omega, w\rangle -\theta(w)\langle \omega, v\rangle - \langle \omega,[v,w]+[v,w]-[v,w]\rangle
=\langle d\omega,\check\xi\wedge\check\eta\rangle.
\end{gather*}

$(ii)$ By applying Proposition~\ref{p31}$(iii)$, we obtain
\begin{gather*}
\langle \frac 1 2[\chi,\chi],\check\xi\wedge\check\eta\rangle =[v,w]-i(\firstbra{v}{\check\eta}_\infty-\firstbra{w}{\check\xi}_\infty
-\rho_1\firstbra{\check\xi}{\check\eta}_\infty)\chi\\
\phantom{\langle \frac 1 2[\chi,\chi],\check\xi\wedge\check\eta\rangle }{}=[v,w]-([v,w]-[w,v]-[v,w])=0.
\end{gather*}

$(iii)$ It follows from  \eqref{buvmod} that, for $\textbf u=\omega\otimes \xi$,
\[
[\chi,\textbf u]=\theta(\chi)\omega\otimes\xi+(-1)^r\omega\wedge [\chi,\xi]-(-1)^{2r}d\omega\wedge i(\xi)\chi
=d\omega\otimes\xi+(-1)^r\omega\wedge[\chi,\xi].
\]
As $D$ is characterized by Proposition~\ref{propriedadesD}, it is enough to prove $[\chi,\xi]=D\xi$. It follows from Propositions~\ref{colchetek}$(ii)$ and \ref{p31}$(ii)$ that
\begin{gather*}
i(\check\eta)[\chi,\xi]=\firstbra{i(\check\eta)\chi}{\xi}_\infty
-i(\firstbra{\check\eta}{\xi}_\infty)\chi=\firstbra{w}{\xi}_\infty=i(\check\eta)D\xi.\tag*{\qed}
\end{gather*}
  \renewcommand{\qed}{}
\end{proof}

If $\textbf u,\textbf v\in \wedge\mathcal T^*\otimes J^\infty\mathcal T$, with $\deg \textbf u=r$, $\deg \textbf v=s$, then we get from \eqref{brauv} and Proposition~\ref{p42}$(iii)$ that
\[D[\textbf u,\textbf v]=[D\textbf u,\textbf v]+(-1)^r[\textbf u,D\textbf v],\]
and
\[[\chi,[\chi,\textbf u]]=[[\chi,\chi],\textbf u]-[\chi,[\chi,\textbf u]]=-[\chi,[\chi,\textbf u]].\]
Therefore, $D^2\textbf u=0$, or
\[D^2=0.\]
Then  \emph{the first linear Spencer complex},
\[0\rightarrow \mathcal T\stackrel{j^\infty}{\rightarrow}J^\infty\mathcal T\stackrel{D}{\rightarrow}\mathcal T^*\otimes J^\infty\mathcal T\stackrel{D}{\rightarrow}\wedge^2\mathcal T^*\otimes J^\infty\mathcal T\stackrel{D}{\rightarrow}\cdots\stackrel{D}{\rightarrow}\wedge^m\mathcal T^*\otimes J^\infty\mathcal T\rightarrow 0,\]
where $\dim T=m$,  is well def\/ined. This complex projects on
\[0\rightarrow \mathcal T\stackrel{j^k}{\rightarrow}J^k\mathcal T\stackrel{D}{\rightarrow}\mathcal T^*\otimes J^{k-1}\mathcal T\stackrel{D}{\rightarrow}\wedge^2\mathcal T^*\otimes J^{k-2}\mathcal T\stackrel{D}{\rightarrow}\cdots\stackrel{D}{\rightarrow}\wedge^m\mathcal T^*\otimes J^{k-m}\mathcal T\rightarrow 0,\]
and is exact (see \cite{ {Ma1},{Ma2},{KS}}).

Let be $\gamma^k$ the kernel of $\pi_k:J^kT\rightarrow J^{k-1}T$. Denote by $\delta$ the restriction of $D$ to $\gamma^k$.   It follows from Proposition~\ref{propriedadesD}$(ii)$ that $\delta$ is $\mathcal{O}_M$-linear and $\delta:\gamma^k\rightarrow T^*\otimes \gamma^{k-1}$.
This map is injective, in fact, if  $\xi\in\gamma^k$, then by \eqref{kernelpi}, $\delta\xi=-\lambda^1(\xi)$ is injective.
As
\[i(v)D(i(w)D\xi)-i(w)D(i(v)D\xi)-i([v,w])D\pi_{k-1}\xi=0,\]
for $v,w\in\mathcal T$, $\xi\in \gamma^k\subset J^k\mathcal T$, we obtain that $\delta$ is symmetric, $i(v)\delta(i(w)\delta\xi)=i(w)\delta(i(v)\delta\xi).$
Observe that we get the map
\[\iota: \ \gamma^k\rightarrow S^2T^*\otimes\gamma^{k-2}\]
def\/ined by
$i(v,w)\iota(\xi)=i(w)\delta(i(v)\delta\xi),$ and, if we go on, we obtain the isomorphism
\[\gamma^k\cong S^kT^*\otimes J^0T,\]
where, given  a basis
$e_1,\dots,e_m\in T$
with  the dual basis $e^1,\dots,e^m\in T^*$,
we obtain the basis
\begin{equation}\label{f}
f^{k_1,k_2,\dots,k_m}_l=\frac{1}{k_1!k_2!\cdots k_m!}(e^1)^{k_1}(e^2)^{k_2}\cdots (e^m)^{k_m}\otimes j^0e_l
\end{equation}
of $S^kT^*\otimes J^0T$, where $k_1+k_2+\cdots+k_m=k$,  $k_1,\dots,k_m\geq 0$ and $l=1,\dots,m$. In this basis
\[\delta\big(f^{k_1,k_2,\dots,k_m}_l\big)=-\sum_{i=1}^me^i\otimes f^{k_1,\dots,k_{i-1},k_i-1,k_{i+1},\dots,k_m}_l.\]
From the linear Spencer complex, we obtain the exact sequence of morphisms of vector bundles
\begin{equation}\label{delta}0\rightarrow \gamma^k\stackrel{\delta}{\rightarrow} T^*\otimes \gamma^{k-1}\stackrel{\delta}{\rightarrow}\wedge^2\mathcal T^*\otimes \gamma^{k-2} \stackrel{\delta}{\rightarrow}\cdots\stackrel{\delta}{\rightarrow}\wedge^m T^*\otimes \gamma^{k-m} \rightarrow 0.\end{equation}

Let's now introduce the \emph{first non-linear Spencer operator} $\mathcal D$. The ``f\/inite'' form $\mathcal D$ of the linear Spencer  operator $D$ is def\/ined by
\begin{equation}\label{dxi}
\mathcal D\sigma=\chi-\sigma_*^{-1}(\chi),
\end{equation}
where $\sigma\in\mathcal Q^{\infty}$.
\begin{proposition}\label{p43}
The operator $\mathcal D$ take values in $\mathcal T^*\otimes J^\infty\mathcal T$, so
\[\mathcal D: \ \mathcal Q^{\infty}\rightarrow \mathcal T^*\otimes J^\infty\mathcal T,\]
and
\begin{equation}\label{ivds}
i(v)(\mathcal D\sigma)_k=\lambda^1\sigma_{k+1}^{-1}.j^1\sigma_k.v-v,
\end{equation}
where $\sigma  = \lim \emph{proj }\sigma_k\in\mathcal Q^{\infty}$.
\end{proposition}
\begin{proof}  By applying \eqref{sigmaomega} and \eqref{sigmastar}, it follows that for  $\xi\in J^\infty\mathcal T$,
\[i(\xi)\mathcal D\sigma=i(\xi)\chi-\sigma^{-1}_*(i(\sigma_*(\xi))\chi)=0,\]
and for $v\in\mathcal T$,
\[i(v)\mathcal D\sigma=i(v)\chi-\sigma^{-1}_*(i(\sigma_*(v))\chi)=v-\sigma_*^{-1}(f_*v),\]
where $f=\beta\circ\sigma$. By Proposition~\ref{actionsigma}$(iii)$,
\begin{equation}\label{ivdsigma}
i(v)(\mathcal D\sigma)_k=v-\big(f_*^{-1}(f_*v)+j^1\sigma_k^{-1}.f_*v.\lambda^1\sigma_{k+1}
-\lambda^1\sigma_{k+1}^{-1}.f_*v.\lambda^1\sigma_{k+1}\big).
\end{equation}
By posing $v=\frac{d}{dt}x_t|_{t=0}$, we obtain
\begin{equation}\label{f41}
j^1\sigma_k.v.j^1\sigma_k^{-1}=\frac{d}{dt}
\big(\sigma_k(x_t).\sigma_k^{-1}(f(x_t))\big)\big|_{t=0}=\frac{d}{dt}f(x_t)\big|_{t=0}=f_*v,
\end{equation}
and by replacing \eqref{f41} in \eqref{ivdsigma}, we get
\begin{gather*}
i(v)(\mathcal D\sigma)_k=-j^1\sigma_k^{-1}.\big(j^1\sigma_k.v.j^1\sigma_k^{-1}\big).\lambda^1\sigma_{k+1}
+\lambda^1\sigma_{k+1}^{-1}.\big(j^1\sigma_k.v.j^1\sigma_k^{-1}\big).\lambda^1\sigma_{k+1} \\
\phantom{i(v)(\mathcal  D\sigma)_k}{}= \big(\lambda^1\sigma_{k+1}^{-1}.j^1\sigma_k.v-v\big).j^1\sigma_k^{-1}.\lambda^1\sigma_{k+1}
 =\big(\lambda^1\sigma_{k+1}^{-1}.j^1\sigma_k.v-v\big).\sigma_k^{-1}.\sigma_{k} \\
\phantom{i(v)(\mathcal  D\sigma)_k}{}
=\lambda^1\sigma_{k+1}^{-1}.j^1\sigma_k.v-v,
\end{gather*}
since that $\lambda^1\sigma_{k+1}^{-1}.j^1\sigma_k.v-v$ is  $\alpha$-vertical (cf.~\eqref{actionvert}).
\end{proof}

\begin{corollary}
We have ${\mathcal D}\sigma=0$ if and only if $\sigma=j^{\infty}(\beta\sigma)$, where $\beta:\mathcal Q^{\infty} \rightarrow \mbox{\rm Dif\/f }M$.
\end{corollary}

\begin{corollary}\label{sigmaD}
If $\sigma_{k+1}\in\mathcal{Q}^{k+1}$, then
\[(\sigma_{k+1})_*(v)=f_*v+(\sigma_{k+1})_*(i(v)\mathcal{D}\sigma_{k+1}),\]
for $v\in \mathcal{T}$.
\end{corollary}
\begin{proof} It follows from \eqref{GTG} and Proposition~\ref{p43} that
\begin{gather*}
(\sigma_{k+1})_*(i(v)\mathcal{D}\sigma_{k+1})=\lambda^1\sigma_{k+1}.(i(v)\mathcal{D}\sigma_{k+1}).\lambda^1
\sigma_{k+1}^{-1}=\lambda^1\sigma_{k+1}.\big(\lambda^1\sigma_{k+1}^{-1}.j^1\sigma_k.v-v\big).\lambda^1\sigma_{k+1}^{-1}\\
\phantom{(\sigma_{k+1})_*(i(v)\mathcal{D}\sigma_{k+1})}{}= j^1\sigma_k.v.\lambda^1\sigma_{k+1}^{-1}
-\lambda^1\sigma_{k+1}.v.\lambda^1\sigma_{k+1}^{-1}=(\sigma_{k+1})_*(v)-f_*v.\tag*{\qed}
\end{gather*}
\renewcommand{\qed}{}
\end{proof}

\begin{proposition}\label{p44}
The operator $\mathcal D$ has the following properties:
\begin{enumerate}\itemsep=0pt
\item [$(i)$] If $\sigma,\sigma'\in \mathcal Q^{\infty}$,
\[\mathcal D(\sigma'\circ\sigma)=\mathcal D\sigma+\sigma_*^{-1}(\mathcal D\sigma').\]
In particular,
\[\mathcal D \sigma^{-1}=-\sigma_*(\mathcal D\sigma).\]
\item [$(ii)$] If $\sigma\in\mathcal Q^{\infty}$, $\mbox{\bf u}\in\wedge\mathcal T^*\otimes J^\infty\mathcal T$,
\[ \mathcal D(\sigma^{-1}_*\mbox{\bf u})=\sigma_*^{-1}(D\mbox{\bf u})+\big[\mathcal D\sigma,\sigma^{-1}_*\mbox{\bf u}\big].\]
\item [$(iii)$] If $\xi=\frac{d}{dt}\sigma_t|_{t=0}$, with $\xi\in J^
\infty\mathcal T$, and $\sigma_t\in\mathcal Q^{\infty}$ is the $1$-parameter group associated to $\xi$, then
\[D\xi=\frac{d}{dt}\mathcal D\sigma_t\big|_{t=0}.\]
\end{enumerate}
\end{proposition}
\begin{proof}
\begin{gather*}
(i)  \ \ \ \mathcal D(\sigma'\circ\sigma)=\chi-\sigma_*^{-1}(\chi)+\sigma_*^{-1}(\chi-(\sigma')_*^{-1}(\chi))=\mathcal D\sigma+\sigma_*^{-1}(\mathcal D\sigma'). \\
(ii)  \ \
 D(\sigma_*^{-1}\textbf u)=\sigma_*^{-1}[\sigma_*\chi,\textbf u]=\sigma_*^{-1}[\chi-\mathcal D\sigma^{-1},\textbf u]
=\sigma_*^{-1}(D\textbf u)-[\sigma_*^{-1}(\mathcal D\sigma^{-1}),\sigma_*^{-1}\textbf u]\\
\phantom{(ii)  \ \  D(\sigma_*^{-1}\textbf u)}{}=\sigma_*^{-1}(D\textbf u)+[\mathcal D\sigma,\sigma_*^{-1}\textbf u].
\\
(iii) \ \frac{d}{dt}\mathcal D\sigma_t\big|_{t=0}=-\frac{d}{dt}(\sigma_t^{-1})_*(\chi)=-[\xi,\chi]=D\xi.\tag*{\qed}
\end{gather*}
\renewcommand{\qed}{}
\end{proof}

Proposition~\ref{p43} says that $\mathcal D$ is projectable:
\begin{gather*}
\mathcal D: \ \mathcal Q^{k+1}\rightarrow \mathcal T^*\otimes J^k\mathcal T,\nonumber\\
\phantom{\mathcal D:}{} \ \ \sigma_{k+1}\mapsto \mathcal D\sigma_{k+1},
\end{gather*}
where \[
i(v)\mathcal D\sigma_{k+1}=\lambda^1\sigma_{k+1}^{-1}.j^1\sigma_k.v-v.\]

It follows from $[\chi,\chi]=0$ that
\[0=\sigma_*^{-1}([\chi,\chi])=[\sigma_*^{-1}(\chi),\sigma_*^{-1}(\chi)]=[\chi-\mathcal D\sigma,\chi-\mathcal D\sigma]=[\mathcal D\sigma,\mathcal D\sigma]-2D(\mathcal D\sigma),\]
therefore
\begin{equation}\label {f43}
D(\mathcal D\sigma)-\frac{1}{2}[\mathcal D\sigma,\mathcal D\sigma]=0.
\end{equation}
If we def\/ine the non-linear operator
\begin{gather*}
\mathcal D_1: \ \mathcal T^*\otimes J^\infty\mathcal T\rightarrow \wedge^2\mathcal T^*\otimes J^\infty\mathcal T, \nonumber\\
\phantom{\mathcal D_1:}{} \ \ \textbf u\mapsto  D\textbf u-\frac{1}{2}[\textbf u,\textbf u],
\end{gather*}
then we can write  \eqref{f43} as
\[\mathcal D_1\mathcal D=0.\]
The operator $\mathcal D_1$ projects in order $k$ to
\[\mathcal D_1: \ \mathcal T^*\otimes J^k\mathcal T\rightarrow \wedge^2\mathcal T^*\otimes J^{k-1}\mathcal T,\]
where
\[\mathcal D_1\textbf u=D\textbf u-\frac{1}{2}\firstbra{\textbf u}{\textbf u}_k.\]
Here $\firstbra{\textbf u}{\textbf u}_k$ denotes  the analogous of formulas \eqref{buv} and \eqref{buvmod} projected in the order $k$, so that the extension of f\/irst bracket makes sense.  We will leave the details to the reader.

We def\/ine the \emph{first non-linear Spencer complex} by
\[1\rightarrow \mbox{Dif\/f }M\stackrel{j^{k+1}}{\rightarrow}\mathcal Q^{k+1}\stackrel{\mathcal D}{\rightarrow} \mathcal T^*\otimes J^k\mathcal T\stackrel{\mathcal D_1}{\rightarrow}\wedge^2\mathcal T^*\otimes J^{k-1}\mathcal T,\]
which is exact in $\mathcal Q^{k+1}$.

It is possible  to def\/ine the f\/irst nonlinear Spencer complex $\mathcal D$ for invertible sections of $Q^{\infty}(M,M')$ by:
\begin{equation}\label{dxig}
\mathcal D\sigma=\chi-\sigma_*^{-1}(\chi'),
\end{equation}
where $\sigma\in\mathcal Q^{\infty}(M,M')$ and $\chi'\in (\check J^\infty\mathcal T')^*\otimes (\check J^\infty\mathcal T')$ is the fundamental form.
The operator $\mathcal D$ take values in $\mathcal T^*\otimes J^\infty\mathcal T$, so
\[\mathcal D: \ \mathcal Q^{\infty}(M,M')\rightarrow \mathcal T^*\otimes J^\infty\mathcal T,\]
and the same formula of Proposition~\ref{p43} holds:
\[i(v)(\mathcal D\sigma)_k=\lambda^1\sigma_{k+1}^{-1}.j^1\sigma_k.v-v,\]
where $\sigma  = \lim \mbox{proj }\sigma_k\in\mathcal Q^{\infty}(M,M')$. Other properties can easily be generalised.

\section{Partial connections}\label{partcon}

In this section we will develop the concept of partial connections or partial covariant derivatives associated with the vector bundle $H\oplus J^kV$ in the directions of the distribution $V\subset T$. We thank the referee for pointing out  that this concept is already in \cite[p.~24]{KT}. The construction of connections for $J^kT$, the transitive case, is in~\cite{G4}.

Let be $V$ an involutive subvector bundle of~$T$,  $a\in M$, and $N$ a (local) submanifold of $M$ such that $T_aN\oplus V_a=T_aM$. Then there exists a coordinate system $(x,y)$ in a neighborhood of $a$ such that $a=(0,0)$, the submanifolds given by points with coordinates $x$ constant are integral submanifolds of $V$, and $N$ is given by the submanifold  $y=0$. At least locally, we can suppose that the integral manifolds of $V$ are the f\/ibers of a submersion $\rho:M\rightarrow N$. In the coordinates~$(x,y)$, we get $\rho(x,y)=(x,0)$.
If we denote by $H$ the subvector bundle of $T$ given by vectors tangent to the submanifolds def\/ined by $y$ constant, then $H$ is involutive and $T=H\oplus V$.  Also, $T^*=	H^*\oplus V^*$. We denote by $\mathcal H$ and $\mathcal V$ the sheaves of germs of~$H$ and~$V$, respectively.

We denote by $Q^k_V$ the subgroupoid of $Q^k$ whose elements are the $k$-jets of local dif\/feomor\-phisms~$h$ of $M$ which are $\rho$-projectable on the identity of $N$. In the coordinates $(x,y)$, $h(x,y)=(x,h_2(x,y))$. The sheaf of germs of  invertible local $\alpha$-sections of $Q^k_V$ will be denoted by $\mathcal{Q}^k_V$. The algebroid 
associated with $Q^k_V$ is $J^kV$, and we denote by $J^k\mathcal V$ the  sheaf of germs  of local sections of $J^kV$. Then the f\/irst non-linear Spencer operator $\mathcal D$ can be restricted to $\mathcal{Q}^{k+1}_V$,
 \begin{equation}\label{first}
 \mathcal{D}: \ \mathcal{Q}^{k+1}_V\rightarrow \mathcal{T}^*\otimes {J}^k\mathcal V,\end{equation}
 and the linear Spencer operator $D$ can be restricted to $J^{k+1} \mathcal V$,
 \begin{equation*}
 {D}: \ J^{k+1}\mathcal V\rightarrow \mathcal{T}^*\otimes {J}^k\mathcal V.
 \end{equation*}

  A vector $u\in T$ decomposes in $u=u_H+u_V$, $u_H\in H$, $u_V\in V$. If $d$ is the exterior dif\/ferential, we get the decomposition $d=d_H+d_V$.
The fundamental form
 $\chi$ decomposes in $\chi=\chi_H+\chi_V$, where $\chi_H(u)=\chi(u_H)$ and $\chi_V(u)=\chi(u_V)$.
  The linear Spencer operator $D$  also decomposes in $D=D_H+D_V$, and it follows from  Proposition~\ref{p42}$(iii)$ that
\[D_H\xi=[\chi_H,\xi]\qquad \mbox{and}\qquad D_V\xi=[\chi_V,\xi],\] for $\xi\in J^k\mathcal{T}$.

 \begin{proposition}\label{p42HV}
If $\omega\in\wedge\mathcal V^*$, and $\mbox{\bf u}\in\wedge\mathcal V^*\otimes J^\infty\mathcal T$, then:
\begin{enumerate}\itemsep=0pt
\item [$(i)$] 
$\theta(\chi_V)\omega=d_V\omega$;
\item [$(ii)$] $[\chi_V,\chi_V]=0$;
\item [$(iii)$] 
$[\chi_V,\mbox{\bf u}]=D_V\mbox{\bf u}$.
\end{enumerate}
\end{proposition}
\begin{proof} Let be $\check\xi=v+\xi,\check\eta=w+\eta\in\mathcal T\oplus J^\infty\mathcal T$, $v=v_H+v_V,w=w_H+w_V\in \mathcal H\oplus \mathcal V$.

$(i)$ As $\theta(\chi_V)$ is a derivation of degree 1, it is enough to prove $(i)$ for $0$-forms $f$ and $1$-forms $\omega\in\mathcal V^*$. From \eqref{theta(u)} we have
\[\theta(\chi_V)f=i(\chi_V)df=d_Vf.\]
 It follows from Proposition~\ref{p31}$(i)$ that
\begin{gather*}
\langle \theta(\chi_V)\omega,\check\xi\wedge\check\eta\rangle=
\theta(v_V)\langle \omega,\check\eta\rangle-\theta(w_V)\langle \omega,\check\xi\rangle
 -\langle \omega,\firstbra{v_V}{\check\eta}_\infty+\firstbra{\check\xi}{w_V}_\infty
-\chi_V(\firstbra{\check\xi}{\check\eta}_\infty)\rangle\\
\phantom{\langle \theta(\chi_V)\omega,\check\xi\wedge\check\eta\rangle} {}
=\theta(v_V)\langle \omega, w_V\rangle-\theta(w_V)\langle \omega, v_V\rangle- \langle \omega,[v_V,w]_V+[v,w_V]_V-[v,w]_V\rangle\\
\phantom{\langle \theta(\chi_V)\omega,\check\xi\wedge\check\eta\rangle} {}
=\langle d\omega,v_V\wedge w_V\rangle =\langle d_V\omega,\check\xi\wedge\check\eta\rangle,
\end{gather*}
since that
\begin{gather*}
[v_V,w]_V+[v,w_V]_V-[v,w]_V\\
\qquad{}=[v_V,w_V]+[v_V,w_H]_V+[v_H,w_V]_V+[v_V,w_V]-[v_H,w_V]_V-[v_V,w_H]_V-[v_V,w_V]\\
\qquad{}=[v_V,w_V].
\end{gather*}

$(ii)$ By applying Proposition~\ref{p31}$(iii)$, we obtain
\begin{gather*}
\langle \frac 1 2[\chi_V,\chi_V],\check\xi\wedge\check\eta\rangle=
[v_V,w_V]-i(\firstbra{v_V}{\check\eta}_\infty-\firstbra{w_V}{\check\xi}_\infty-
(\rho_1)_*\firstbra{\check\xi}{\check\eta}_\infty)\chi_V \\
\phantom{\langle \frac 1 2[\chi_V,\chi_V],\check\xi\wedge\check\eta\rangle}{}
=[v_V,w_V]-([v_V,w]-[w_V,v]-[v,w])_V=0.
\end{gather*}

$(iii)$ It follows from  \eqref{buvmod} that
\begin{gather*}
[\chi_V,\textbf u]=\theta(\chi_V)\omega\otimes\xi+(-1)^r\omega\wedge [\chi_V,\xi]-(-1)^{2r}d\omega\wedge i(\xi)\chi_V \\
\phantom{[\chi_V,\textbf u]}{}=d_V\omega\otimes\xi+(-1)^r\omega\wedge[\chi_V,\xi].
\end{gather*}
It is enough to prove $[\chi_V,\xi]=D_V\xi$. It follows from Propositions~\ref{colchetek}$(ii)$ and \ref{p31}$(ii)$ that
\begin{gather*}
i(\check\eta)[\chi_V,\xi]=
\firstbra{i(\check\eta)\chi_V}{\xi}_\infty-i(\firstbra{\check\eta}{\xi}_\infty)\chi_V
=\firstbra{w_V}{\xi}_\infty=i(\check\eta)D_V\xi.\tag*{\qed}
\end{gather*}
\renewcommand{\qed}{}
\end{proof}

From item $(ii)$ of this proposition, \eqref{dtutv} and \eqref{tutv} we obtain
\begin{equation}\label{dV2}
d_V^2=0.
\end{equation}
  The f\/irst non-linear Spencer operator $\mathcal{D}$ also decomposes  naturally in $\mathcal{D}=\mathcal{D}_H+\mathcal{D}_V$, where $i(u)\mathcal{D}_H=i(u_H)\mathcal{D}$ and $i(u)\mathcal{D}_V=i(u_V)\mathcal{D}$. We obtain:
\begin{proposition}\label{dec}If $F\in\mathcal{Q}^{k+1}_V$ is such that $f=\beta F$ satisfies $f_*H=H$, then
\[\mathcal{D}_HF=\chi_H-F^{-1}_*(\chi_H),\]
and
\[\mathcal{D}_VF=\chi_V-F^{-1}_*(\chi_V).\]
\end{proposition}

\begin{proof}  From hypothesis, we get $f_*(u_H)=(f_*u)_H$ and $f_*(u_V)=(f_*u)_V$, and from equa\-tion~\eqref{dxi}, we get
\begin{gather*}
i(u)\mathcal{D}_HF=i(u_H)\mathcal{D}F=i(u_H)\chi-F_*^{-1}(i(f_*(u_H))\chi)=
i(u_H)\chi-F_*^{-1}(i((f_*u)_H)\chi)\vspace{.2cm}\\
\phantom{i(u)\mathcal{D}_HF}{}=i(u)\chi_H-F_*^{-1}(i(f_*u)\chi_H)=i(u)(\chi_H-F^{-1}_*(\chi_H)),
\end{gather*}
so $\mathcal{D}_HF=\chi_H-F^{-1}_*(\chi_H)$. The proof of the second formula is analogous. \end{proof}

If we apply Proposition~\ref{dec} to $GF$, we get
\[
\mathcal{D}_V(GF)=\chi_V-F^{-1}_*G^{-1}_*(\chi_V)=\chi_V+F^{-1}_*(\mathcal{D}_VG-\chi_V)=\mathcal{D}_VF+F^{-1}_*(\mathcal{D}_VG),
\]
and if we pose $F=G^{-1}$, we get
\begin{equation*}
\mathcal{D}_VG^{-1}=-G_*(\mathcal{D}_VG).
\end{equation*}


\begin{definition}
A \emph{partial connection} $\nabla$ on $H\oplus J^kV$ is a $\mathbb{R}$-linear map
\[\nabla:\mathcal{H}\oplus J^k\mathcal{V}\rightarrow \mathcal{V}^*\otimes (\mathcal{H}\oplus J^k\mathcal{V})\]
such that
\[\nabla(f\check\xi)=(d_Vf)\otimes\check\xi+f\nabla\xi,\]
for $f\in \mathcal{O}_M$, $\check\xi\in\mathcal{H}\oplus J^k\mathcal{V}$.
We extend $\nabla$ to $\wedge V^*\otimes(H\oplus J^kV)$ by
\begin{equation*}
\nabla(\alpha\otimes\check\xi)=d_V\alpha\otimes\check\xi+(-1)^{| \alpha|}\alpha\wedge\nabla\check\xi,
\end{equation*}
where $\alpha\in \wedge \mathcal{V}^*$ and $|\alpha|$ is the degree of $\alpha$.
\end{definition}
It follows from \eqref{dV2} that $\nabla^2(\alpha\otimes \check\xi)=\alpha\wedge\nabla^2\check\xi$, so
\[\nabla^2:\mathcal{H}\oplus J^k\mathcal{V}\rightarrow \wedge^2\mathcal{V}^*\otimes (\mathcal{H}\oplus J^k\mathcal{V})\] is a tensor, called the \emph{curvature tensor} of the partial connection $\nabla$. If $\nabla^2=0$ we say that $\nabla$ is \emph{flat}.

Let be $\omega\in{\mathcal{V}}^*\otimes J^{k+1}\mathcal V$ such that $\beta_*(i(w)\omega)=w$, for $w\in \mathcal{V}$, and
\[ \tilde \omega=\chi_V+\omega,\]
so $\tilde\omega\in{\mathcal{V}}^*\otimes J^{k+1}\tilde{\mathcal V}$, where $J^{k+1}\tilde{\mathcal V}=J^{k+1}\tilde{\mathcal T}\cap (\mathcal T\oplus J^{k+1}\mathcal V)$.

In the sequel, for $\check\xi\in \mathcal{H}\oplus {J}^k\mathcal{V}$, we denote by $\thirdbra{\tilde\omega}{\check\xi}_k$ and $\secondbra{\tilde\omega}{\tilde\omega}_{k+1}$  the analogous of formu\-las~\eqref{buv} and \eqref{buvmod} projected in the order $k$ and $k+1$, respectively, so that the same construction of third and second bracket make sense. We will use this convention in the present section when it makes sense, and leave the details to the reader.
\begin{proposition}\label{opnabla}The operator $\nabla$ defined by
\[\nabla\check\xi=\thirdbra{\tilde\omega}{\check\xi}_k,\]
for $\tilde\omega\in\mathcal{V}^*\otimes  J^{k+1}\tilde{\mathcal{V}}$, $\check\xi=u+\xi\in \mathcal{H}\oplus {J}^k\mathcal{V}$ is a partial connection on $H\oplus J^kV$ with curvature \[\nabla^2=\frac12\secondbra{\tilde\omega}{\tilde\omega}_{k+1}.\]
\end{proposition}
\begin{proof}  If we apply formula $(ii)$ of Proposition~\ref{p31}, we obtain
\[
i(\check\eta)\thirdbra{\tilde\omega}{\check\xi}_k
=\thirdbra{i(\check\eta)\tilde\omega}{\check\xi}_k
-i(\firstbra{\check\eta_{k+1}}{\check\xi_{k+1}}_{k+1})\pi_k\tilde\omega,\]
where $\check\eta=v+\eta\in\mathcal{T}\oplus {J}^k\mathcal{T}$, and $\check\eta_{k+1},\check\xi_{k+1}\in\mathcal{T}\oplus {J}^{k+1}\mathcal{T}$ projects on $\check\eta$, $\check\xi$, respectively. If $v\in\mathcal{H}$, the right side is $0$, and if $v\in\mathcal{V}$
\[i(\check\eta)\thirdbra{\tilde\omega}{\check\xi}_k=[v,u]_H+\thirdbra{i(v)\tilde\omega}{\xi}_k-i([v,u])\pi_k\omega-i(u)D(i(v)\omega)\in\mathcal{H}\oplus J^k\mathcal{V}.\]
Then $\nabla\check\xi\in\mathcal V^*\otimes J^k\mathcal V$.
Also,
\begin{gather*}
i(\check\eta)\thirdbra{\tilde\omega}{f\check\xi}_k=
\thirdbra{i(\check\eta)\tilde\omega}{f\check\xi}_k-i(\firstbra{\check\eta_{k+1}}{f\check\xi_{k+1}}_{k+1})
\pi_k\tilde\omega \\
\phantom{i(\check\eta)\thirdbra{\tilde\omega}{f\check\xi}_k}{}
=v_V(f)\check\xi+f\thirdbra{i(\check\eta)\tilde\omega}{\check\xi}_k
-i(v(f)\check\xi+f\firstbra{\check\eta_{k+1}}{\check\xi_{k+1}}_{k+1})\pi_k\tilde\omega \\
\phantom{i(\check\eta)\thirdbra{\tilde\omega}{f\check\xi}_k}{}
=v_V(f)\check\xi+f\thirdbra{i(\check\eta)\tilde\omega}{\check\xi}_k
-i(\firstbra{\check\eta_{k+1}}{\check\xi_{k+1}}_{k+1})(f\pi_k\tilde\omega)\\
\phantom{i(\check\eta)\thirdbra{\tilde\omega}{f\check\xi}_k}{}
=v_V(f)\check\xi+i(\check\eta)(f\thirdbra{\tilde\omega}{\check\xi}_k)
 =i(\check\eta)(d_Vf)\check\xi+i(\check\eta)(f\thirdbra{\tilde\omega}{\check\xi}_k),
\end{gather*}
so
\[\nabla(f\check\xi)=d_Vf\otimes\check\xi+f\nabla\check\xi.\]
If $\alpha\otimes\check\xi\in\wedge\mathcal V^*\otimes (\mathcal H\oplus J^k\mathcal V)$, we know  from \eqref{buvmod} that 
\[\thirdbra{\tilde\omega}{\alpha\otimes\check\xi}_k=\theta(\tilde\omega)\alpha\otimes\check\xi+(-1)^{|\alpha|}\alpha\wedge\thirdbra{\tilde\omega}{\check\xi}_k-(-1)^{2|\alpha|}d\alpha\wedge i(\check\xi){\tilde\omega},
\]
and from \eqref{theta(u)} that
\[\theta(\tilde\omega)\alpha=i(\tilde\omega)d\alpha-d(i(\tilde\omega)\alpha)=i(\chi_V)d\alpha-d(i(\chi_V)\alpha)=\theta(\chi_V)\alpha=d_V\alpha,\]
and as $i(\check\xi){\tilde\omega}=0$, then
\[\thirdbra{\tilde\omega}{\alpha\otimes\check\xi}_k=d_V\alpha\otimes\check\xi+(-1)^{|\alpha|}\alpha\wedge\thirdbra{\tilde\omega}{\check\xi}_k=d_V\alpha\otimes\check\xi+(-1)^{|\alpha|}\alpha\wedge\nabla\check\xi,
\]
so
\[\nabla(\alpha\otimes\check\xi)=\thirdbra{\tilde\omega}{\alpha\otimes\check\xi}_k.\]
Therefore, \[\nabla^2\check\xi=\thirdbra{\tilde\omega}{\thirdbra{\tilde\omega}{\check\xi}_k}_k=\thirdbra{\secondbra{\tilde\omega}{\tilde\omega}_{k+1}}{\check\xi}_k-\nabla^2\check\xi,\] and from this it follows that
\begin{gather*}
\nabla^2\check\xi=\frac12\thirdbra{ \secondbra{\tilde\omega}{\tilde\omega}_{k+1}}{\check\xi}_k.\tag*{\qed}
\end{gather*}
\renewcommand{\qed}{}
\end{proof}

Let be $\sigma^y:N\rightarrow Q^{k+1}_V$ a family of dif\/ferentiable sections such that
\begin{equation}\label{sigmay}
\sigma^y(x,0)\in Q^{k+1}_V((x,0),(x,y)),
\end{equation}
with $\sigma^0(x,0)=j^{k+1}_{(x,0)}\mbox{id}$. Then $\sigma:M\rightarrow Q^{k+1}_V$, given by  $\sigma(x,y)=\sigma^y(x,0)$, is a dif\/ferentiable $\beta$-section of $Q^{k+1}_V|_N=\{X\in Q^{k+1}_V:\alpha(X)\in N\}$. If $v(x,y)=\frac{d}{dt}(x,y(t))|_{t=0}$, def\/ine
\[\omega\in\mathcal{V}^*\otimes J^{k+1}\mathcal{V}\]
by
\begin{equation}\label{omegacon}i(v(x,y))\omega
=\frac{d}{dt}\big(\sigma^{y(t)}(x,0)\sigma^y(x,0)^{-1}\big)\big|_{t=0}\in J^{k+1}_{(x,y)}V.\end{equation}
\begin{proposition}\label{flat}
The partial connection $\nabla$ defined by $\tilde\omega=\chi_V+\omega$, where $\omega$  is defined as in~\eqref{omegacon}, is flat.
\end{proposition}

\begin{proof}
First of all,
\[\beta_*(i(v(x,y))\omega)=\frac{d}{dt}\beta\big(\sigma^{y(t)}(x,0)\sigma^y(x,0)^{-1}\big)\big|_{t=0}
=\frac{d}{dt}(x,y(t))|_{t=0}=v(x,y),\]
 so $\omega$ satisf\/ies the condition to def\/ine a partial connection. We will show that $\tilde\omega$ satisf\/ies $\secondbra{\tilde\omega}{\tilde\omega}_{k+1}=0$, i.e., the partial connection $\nabla$ def\/ined by $\tilde\omega$ is f\/lat. Let be, for $v\in \mathcal{V}$, $\overline{i(v)\omega}$ the right invariant vector f\/ield def\/ined by
 \[\overline{i(v)\omega}(X)=i(v)\omega(\beta(X)).X,\]
  where $X\in Q^{k+1}_V$. We prove that $\overline{i(v)\omega}$ is tangent to $\sigma(M)$, which follows from
\begin{gather*}
\overline{i(v)\omega}(\sigma(x,y))=i(v(x,y))\omega.\sigma(x,y)
=\frac{d}{dt}\big[\sigma^{y(t)}(x,0)\sigma^y(x,0)^{-1}\big]_{t=0}.\sigma(x,y) \\
\phantom{\overline{i(v)\omega}(\sigma(x,y))} {}=\frac{d}{dt}\big[\sigma^{y(t)}(x,0)\sigma^y(x,0)^{-1}\sigma(x,y)\big]_{t=0}
=\frac{d}{dt}\big[\sigma^{y(t)}(x,0)\big]_{t=0}\in T(\sigma(M)).
\end{gather*}
To f\/inish, we know that, for $v,w\in\mathcal V$, $\overline{i(v)\omega}$, $\overline{i(w)\omega}$ and  $\overline{i([v,w])\omega}$ are tangent to the submani\-fold~$\sigma(M)$, and, as $\beta_*([\overline{i(v)\omega},\overline{i(w)\omega}]=[v,w]$, it follows that
\[[\overline{i(v)\omega}, \overline{i(w)\omega}]=\overline{i([v,w])\omega},\]
so from Proposition \ref{tilde},
\[\secondbra{i(v)\tilde\omega}{i(w)\tilde\omega}_{k+1}=i([v,w])\tilde\omega.\]
From Proposition~\ref{p31}$(iii)$ we obtain
\begin{gather*}
 i( v\wedge w)\Big(\frac 1 2\secondbra{\tilde\omega}{\tilde\omega}_{k+1}\Big)\\
\qquad{}=\secondbra{i(v)\tilde\omega}{i(w)\tilde\omega}_{k+1}-i\left(\firstbra{i(v)\tilde\omega_{k+2}}{w}_{k+2}
-\firstbra{i(w)\tilde\omega_{k+2}}{v}_{k+2}-i([v,w])\tilde\omega\right)\tilde\omega \\
\qquad{}=\secondbra{i(v)\tilde\omega}{i(w)\tilde\omega}_{k+1}-i\left([v,w]-[w,v]-[v,w]\right)\tilde\omega
=0,
\end{gather*}
where $\omega_{k+2}$ is a section in $J^{k+2}\mathcal V$ that projects on $\omega$. The proposition is proved.
\end{proof}

Therefore, given a section $u+\xi:N\rightarrow(H\oplus J^{k}V)|_N$, there exists only one section $(U+\Xi)\in \mathcal{H}\oplus J^{k}\mathcal{V}$ such that $(U+\Xi)|_N=u+\xi$ and $\nabla(U+\Xi)=0$. The following proposition characterizes these sections:
\begin{proposition}\label{UXI}Let be  $(U+\Xi)\in \mathcal{H}\oplus J^{k}\mathcal{V}$ such that
\[(U+\Xi)|_N=u+\xi:N\rightarrow(H\oplus J^{k}V)|_N,\]
 and $\nabla(U+\Xi)=0$. Then
\begin{equation}\label{UXi}
(U+\Xi)(x,y)=\sigma^y_*((u+\xi)(x,0)).
\end{equation}
\end{proposition}
\begin{proof} Choose $v\in\mathcal{V}$ such that $v$ is $H$-projectable, i.e., the 1-parameter group $f_t$ of $v$ is given in coordinates by $f_t(x,y)=(x,h_t(y))$. If we def\/ine
\[F_t(x,y)=\sigma^{h_t(y)}(x,0)\sigma^y(x,0)^{-1},\]
then $F_t\in Q^{k+1}_V$. Furthermore, $\beta F_t=f_t$, and
\begin{gather*}
(F_sF_t)(x,y)=\big(\sigma^{h_s(h_t(y))}(x,0)\sigma^{h_t(y)}(x,0)^{-1}\big)
\big(\sigma^{h_t(y)}(x,0)\sigma^y(x,0)^{-1}\big) \\
\phantom{(F_sF_t)(x,y)}{}=\sigma^{h_{s+t}(y)}(x,0)\sigma^y(x,0)^{-1} =F_{s+t}(x,y).
\end{gather*}
So, $F_t$ is the 1-parameter group such that $\frac{d}{dt}F_t|_{t=0}=\tilde\omega(v)$. If $U+\Xi$ is def\/ined by \eqref{UXi}, we get
\begin{gather*}
((F_t)_*(U+\Xi))(x,y)=(F_t)_*((U+\Xi)(x,h_{-t}(y)))
=(F_t)_*(\sigma^{h_{-t}(y)}_*((u+\xi)(x,0))) \\
\phantom{((F_t)_*(U+\Xi))(x,y)}{}=\big(\sigma^{h_t(h_{-t}(y))}\big(\sigma^{h_{-t}(y)}\big)^{-1}\big)_*
 \big(\sigma^{h_{-t}(y)}\big)_*((u+\xi)(x,0))
 \\
\phantom{((F_t)_*(U+\Xi))(x,y)}{}=\sigma_*^y((u+\xi)(x,0))=(U+\Xi)(x,y),
\end{gather*}
so $\thirdbra{i(v)\tilde\omega}{U+\Xi}_k=0$. Also,
\[(U+\Xi)(x,0)=\sigma^0_*((u+\xi)(x,0))=(j^{k+1}_{(x,0)}\mbox{id})_*((u+\xi)(x,0))=(u+\xi)(x,0).\]
 Let $\bar u\in\mathcal{H}$ be the vector f\/ield $\rho$-projectable such that $\bar u|_N=u$. Then \[[\bar u,v]=0,\] and from \eqref{UXi} and Corollary~\ref{sigmaD} we get
 \[U(x,y)=\bar u(x,y)+\sigma^y_*\left(i(u(x,0))\mathcal{D}\sigma^y\right)\in\mathcal{H}\oplus J^k\mathcal{V},\] so
\begin{gather*}
i(v)\nabla(U+\Xi)=
\thirdbra{i(v)\tilde\omega}{U+\Xi}_k-i([v,\bar u])\pi_k\tilde\omega=0.\tag*{\qed}
\end{gather*}
\renewcommand{\qed}{}
\end{proof}

We will now verify how  a partial connection def\/ined by $\tilde\omega$ changes.

Let $M'$,  $V'$, $a'$, $N'$,  $(x',y')$, $H'$, $\rho'$, be as above, with the same properties and dimensions. Denote by $T'=TM'$. Let be $\phi:N\rightarrow N'$ a (local) dif\/feomorphism, with $\phi(a)=a'$, and denote by $Q^k_{\phi}$ the submanifold of $Q^k(M,M')$ of $k$-jets of local dif\/feomorphisms $\tau:M\rightarrow M'$ such that $\rho'\tau=\phi\rho$. This means $\tau(x,y)=(\phi(x),b(x,y))$.

Let be $\mathcal{Q}^k_\phi$ the sheaf of germs of invertible local sections of ${Q}^k_\phi$. Then,  by restriction of action~\eqref{GTGduplo},  there exists an action of $\mathcal Q^{k+1}_\phi$ (similar to \eqref{GTG}) on $T\oplus J^k V$:
\begin{gather*}
\mathcal Q^{k+1}_\phi\times(\mathcal T\oplus J^k\mathcal V)\rightarrow \mathcal T'\oplus J^k\mathcal V',\\
 (\sigma_{k+1},v+\xi_k)\mapsto  (\sigma_{k+1})_*(v+\xi_k).
\end{gather*}

The operator $\mathcal D:\mathcal Q^{k+1}(M,M')\rightarrow \mathcal T^*\otimes J^k\mathcal T$ def\/ined in \eqref{dxig} restricts, as  \eqref{first}, to
 \begin{equation*}
 \mathcal{D}:\mathcal{Q}^{k+1}_{\phi}\rightarrow \mathcal{T}^*\otimes {J}^k\mathcal V,
 \end{equation*}
 and as above $\mathcal{D}$ decomposes in $\mathcal{D}=\mathcal{D}_H+\mathcal{D}_V$. The analogous of Proposition \ref{dec} holds:
 \begin{proposition}\label{decphi}
 If $F\in\mathcal{Q}^{k+1}_{\phi}$ is such that $f=\beta F$ satisfies $f_*H=H'$, then
 \[\mathcal{D}_HF=\chi_H-F^{-1}_*(\chi_{H'}),\]
 and
 \[ \mathcal{D}_VF=\chi_V-F^{-1}_*(\chi_{V'}).\]
\end{proposition}
Denote by
\[
\mathcal{Q}^{k+1}_{\phi^{-1}}=\big\{\Phi\in \mathcal Q^{k+1}(M',M):\Phi^{-1}\in\mathcal{Q}^{k+1}_{\phi}\big\}.
\]
If $F\in \mathcal{Q}^{k+1}_{\phi}$ satisfy for $f=\beta F$, $f_*(H)=H'$, then
\begin{equation}\label{dvlinha}
\mathcal{D}_{V'}F^{-1}=\chi_{V'}-F_*(\chi_{V})
\end{equation}
and if $G\in \mathcal{Q}^{k+1}_{\phi^{-1}}$
satisfy $\alpha(G)=\beta(F)$, and for $g=\beta G$ we have $g_*(H')=H$, then by applying Propositions \ref{dec} and  \ref{decphi} to $FG$ we get
\[
\mathcal{D}_{V'}(FG)=\chi_{V'}-G^{-1}_*F^{-1}_*(\chi_{V'})=\chi_{V'}+G^{-1}_*(\mathcal{D}_{V}F-\chi_{V})=\mathcal{D}_{V'}G+G^{-1}_*(\mathcal{D}_{V}F).
\]

By posing $G=F^{-1}$ we get
\begin{equation}\label{DVGphi}
\mathcal{D}_{V'}F^{-1}=-F_*(\mathcal{D}_{V}F).
\end{equation}

Choose $\Phi\in\mathcal{Q}^{k+1}_\phi$ with $\varphi=\pi_0\Phi$ satisfying $\varphi(N)=N'$ and $\varphi_*(H)=H'$. Then  $\varphi|_N=\phi$ and $(x',y')=\varphi(x,y)=(\phi(x),b(y))$. Def\/ine $\sigma'^{y'}$ as
\begin{equation}\label{sigmalinha}
\sigma'^{y'}(x',0)=\Phi(x,y)\sigma^y(x,0)\Phi(x,0)^{-1},
\end{equation}
and let be $\omega'$ and $\nabla'$ as in \eqref{omegacon} and Proposition~\ref{opnabla}, respectively. Following the proof of Proposition~\ref{UXI}, take $v\in\mathcal{V}$ such that $v$ is $H$-projectable, i.e., the 1-parameter group $f_t$ of $v$ is given in coordinates by $f_t(x,y)=(x,h_t(y))$. Def\/ine $F_t\in Q^{k+1}_V$ by
\[F_t(x,y)=\sigma^{h_t(y)}(x,0)\sigma^y(x,0)^{-1}.\]If $v'=\varphi_*v$, then $v'\in \mathcal{V}'$ is $H'$-projectable, $f'_t=\varphi f_t \varphi^{-1}$ is the associated 1-parameter group of~$v'$, and the 1-parameter group associated with $\omega'(v')$ satisf\/ies
\begin{gather*}
F'_{t'}(x',y')=\sigma'^{h_{t'}(y')}(x',0)\sigma^{y'}(x',0)^{-1} \\
\phantom{F'_{t'}(x',y')}{}=\big(\Phi(x,h_t(y))\sigma^{h_t(y)}(x,0)\Phi(x,0)^{-1}\big)
\big(\Phi(x,y)\sigma^y(x,0)\Phi(x,0)^{-1}\big)^{-1} \\
\phantom{F'_{t'}(x',y')}{}=\Phi(x,h_t(y))F_t(x,y)\Phi(x,y)^{-1}  =(\Phi F_t\Phi^{-1})(x',y'),
\end{gather*}
i.e., $F'_t=\Phi F_t\Phi^{-1}$. From this, we get $i(v')\tilde\omega'=\Phi_*(i(v)\tilde\omega)$, and
\[i(v')(\Phi_*\tilde\omega)=\Phi_*(i(\Phi_*^{-1}(v'))\tilde\omega)
=\Phi_*(i(f_*^{-1}(v'))\tilde\omega)=\Phi_*(i(v)\tilde\omega)=i(v')\tilde\omega',\]
so
\begin{equation}
\label{phiom}\Phi_*\tilde\omega=\tilde\omega'.
\end{equation}
Then
\[\nabla'(\Phi_*(U+\Xi))=\thirdbra{\tilde\omega'}{\Phi_*(U+\Xi)}_k=\Phi_*\thirdbra{\tilde\omega}{U+\Xi}_k,\]
i.e.,
\begin{equation}\label{nabla1}
\nabla'\Phi_*=\Phi_*\nabla,
\end{equation}
which shows  that $\Phi_*(U+\Xi)$ is parallel with respect to $\nabla'$ if  and only if $U+\Xi$ is parallel with respect to $\nabla$.

Taking account of \eqref{dvlinha}, the equation \eqref{phiom} projected  in order $k$ is equivalent to,
\[\pi_k(\chi_{V'}+\omega')=\pi_k(\Phi_*(\chi_V+\omega))=\Phi_*(\chi_V)+\Phi_*(\pi_k\omega)=
\big(\chi_{V'}-\mathcal{D}_{V'}\Phi^{-1}\big)+\Phi_*(\pi_k\omega),\]
or, considering \eqref{DVGphi},
\[
\chi_{V'}+\pi_k\omega'=\chi_{V'}+\Phi_*(\mathcal{D}_V\Phi)+\Phi_*(\pi_k\omega),
\]
so
\begin{equation}\label{pik}
\pi_k\omega'=\Phi_*(\mathcal{D}_V\Phi+\pi_k\omega).
\end{equation}

\section{Linear Lie equations}\label{sLieeq}
\begin{definition}\label{prolongation}
Let be $R^k$ a subvector bundle of $J^kT$. We def\/ine the \emph{prolongation} $R^{k+1}$ of $R^k$ by
\[R^{k+1}=(\lambda^1)^{-1}(J^1R^k\cap \lambda^1(J^{k+1}T))\subset J^{k+1}T,\]
where the intersection is done in $J^1J^kT$.
\end{definition}
We denote the prolongation $(R^{k+1})^{+1}$ of $R^{k+1}$ by $R^{k+2}$ and so on, and  by $\mathcal R^{k+l}$ the sheaf of germs  of local sections of $R^{k+l}$, for $l\geq 0$.
\begin{proposition} \label{prolongationx} A section $\xi\in J^{k+1}\mathcal T$ is in $\mathcal R^{k+1}$ if and only if $\pi_k\xi\in \mathcal R^k$ and $D\xi\in\mathcal T^*\otimes \mathcal R^k$.
\end{proposition}
\begin{proof}
From Def\/inition~\ref{prolongation} and \eqref{D}, since  $j^1\pi_k\xi\in J^1\mathcal{R}^k$.
\end{proof}

\begin{definition}\label{Lieequation}
A subvector bundle $R^k$ of $J^kT$ is a \emph{linear Lie equation} if the prolongation $R^{k+1}$ of $R^k$ is  a subvector bundle of $J^{k+1}T$ such that
\begin{enumerate}\itemsep=0pt
\item[$(i)$]$\pi_k(R^{k+1})=R^k$;
\item[$(ii)$]$\firstbra{\mathcal R^{k+1}}{\mathcal R^{k+1}}_{k+1}\subset{\mathcal R^{k}}.$
\end{enumerate}
\end{definition}
It follows from Proposition~\ref{prolongationx} and Def\/inition~\ref{Lieequation} that
\begin{equation*}
\firstbra{\mathcal T\oplus \mathcal R^{k+1}}{\mathcal T\oplus \mathcal R^{k+1}}_{k+1}\subset\mathcal T\oplus \mathcal R^{k},
\end{equation*}
and from this,
\begin{equation}
\label{lietilde}\secondbra{\mathcal{\tilde R}^k}{\mathcal {\tilde R}^k}_k\subset\mathcal {\tilde R}^k,
\end{equation}
and
\begin{equation}\label{con}\thirdbra{\mathcal {\tilde R}^{k+1}}{\mathcal R^k}_k \subset\mathcal R^k.\end{equation}

\begin{proposition}
If $R^k\subset J^kT$ is a linear Lie equation, then
\[\firstbra{\mathcal T\oplus \mathcal R^{k+l}}{\mathcal T\oplus \mathcal R^{k+l}}_{k+l}\subset {\mathcal T\oplus \mathcal R^{k+l-1}},\]
for $l\geq 2$.
\end{proposition}
\begin{proof} Let's prove this for $l=2$.  The other proofs for $l>2$ are equal. Suppose $\xi_{k+2},\eta_{k+2}\in\mathcal R^{k+2}$. Then, by Proposition~\ref{prolongationx},  $\xi_{k+1}\!=\pi_{k+1}\xi_{k+2},\eta_{k+1}\!=\pi_{k+1}\eta_{k+2}\in\mathcal R^{k+1}$, and $D\xi_{k+2},D\eta_{k+2}\in\mathcal T\otimes \mathcal R^{k+1}$. So,
\[\pi_k\firstbra{\xi_{k+2}}{\eta_{k+2}}_{k+2}=\firstbra{\xi_{k+1}}{\eta_{k+1}}_{k+1}\in\mathcal R^k,\]
and
\[D\firstbra{\xi_{k+2}}{\eta_{k+2}}_{k+2}=\firstbra{D\xi_{k+2}}{\eta_{k+1}}_{k+1}+\firstbra{\xi_{k+1}}{D\eta_{k+2}}_{k+1}\in\mathcal T^*\otimes \mathcal R^k.\]
 Therefore, by the same Proposition~\ref{prolongationx},
$\firstbra{\xi_{k+2}}{\eta_{k+2}}_{k+2}\in \mathcal R^{k+1}$, and the proposition follows.\end{proof}

It does not  follow from this proposition that $R^{k+l}$ is a vector bundle, and that $\pi_{k+l}:R^{k+l}\rightarrow R^{k+l-1}$ is onto, for $l\geq 2$. To obtain this, we need an additional condition.
\begin{definition} We say that the linear Lie equation $R^k$ is \emph{formally integrable} if
\begin{enumerate}\itemsep=0pt
	\item[$(i)$] $R^{k+l}$ is a subvector bundle of $J^{k+l}T$,
	\item[$(ii)$] $\pi_{k+l}:R^{k+l+1}\rightarrow R^{k+l}$ is onto,
\end{enumerate}
for $l\geq 1$.
\end{definition}
The \emph{symbol} $g^k$ of $R^k$ is the kernel of $\pi_{k-1}:R^k\rightarrow J^{k-1}T$. Also, $g^{k+l}$ is the kernel of $\pi_{k+l-1}:R^{k+l}\rightarrow R^{k+l-1}$, for $l\geq 1$. It follows from Proposition~\ref{prolongationx} and from \eqref{delta} that we have the subcomplex
 \begin{equation}\label{deltas}0\rightarrow g^{k+l}\stackrel{\delta}{\rightarrow} T^*\otimes g^{k+l-1}\stackrel{\delta}{\rightarrow}\wedge^2\mathcal T^*\otimes g^{k+l-2} \stackrel{\delta}{\rightarrow}\wedge^3 T^*\otimes \gamma^{k+l-3} \end{equation}
for $l\geq 2$.
\begin{definition}We say that the symbol $g^k$ is 2-\emph{acyclic} if the subcomplex \eqref{deltas} is exact for $l\geq 2$.
\end{definition}
The following proposition is in \cite{{G1},{G2}}. For an alternative proof, see \cite{{Ru1},{Ru2},{Ru3},{V2}}.
\begin{proposition}
If $R^k\subset J^kT$ is such that
\begin{enumerate}\itemsep=0pt
	\item[$(i)$] $R^{k+1}$ is a subvector bundle of $J^{k+1}T$,
	\item[$(ii)$] $\pi_{k}:R^{k+1}\rightarrow R^{k}$ is onto,
	\item[$(iii)$] $g^k$ is $2$-acyclic,
\end{enumerate}
then $\pi_{k+l-1}:R^{k+l}\rightarrow R^{k+l-1}$ is onto for $l\geq 2$.
\end{proposition}
A consequence of this proposition is:
\begin{corollary}\label{rkfi}If $R^k\subset J^kT$ is a linear Lie equation and $g^k$ is $2$-acyclic, then $R^k$ is formally integrable.
\end{corollary}

Given a linear Lie equation $R^k$, let be the distribution $B\subset TQ^k$ def\/ined by $B_X=R^k_{\beta(X)}.X$, for $X\in Q^k$. It follows from \eqref{Liebra} and \eqref{lietilde} that the distribution $B$ is involutive.  Let be~$P^k(x)$ the integral leaf of $B$ that contains the point $I(x)$, and $P^k=\cup_{x\in M}P^k(x)$. Then $P^k$ is a groupoid, and a dif\/ferentiable submanifold at a neighborhood of $I$. As our problem is local, we will suppose that $P^k$ is a dif\/ferentiable groupoid,   the  dif\/ferentiable groupoid associated with the linear Lie equation~$R^k$. Then the linear Lie equation $R^k$ is the Lie algebroid associated with~$P^k$. As before, we denote by~$\mathcal P^k$ the groupoid of invertible sections of~$P^k$.

We def\/ine the \emph{prolongation} $P^{k+1}$ of $P^k$ by
\[P^{k+1}=(\lambda^1)^{-1}\big(Q^1P^k\cap\lambda^1Q^{k+1}\big),\]
where  $\lambda^1:Q^{k+1}\rightarrow Q^1Q^k$ and $Q^1P^k$ is the groupoid of 1-jets of invertible sections of $P^k$. The following  is Proposition~6.9$(ii)$  of~\cite{Ma2}:

\begin{proposition}\label{prolqk}
Let be $F\in\mathcal Q^{k+1}$ such that
\begin{enumerate}\itemsep=0pt
	\item [$(i)$]$\pi_{k}F\in \mathcal P^k$,
	\item [$(ii)$]$\mathcal DF\in\mathcal T^*\otimes \mathcal R^{k}$.
\end{enumerate}
Then $F\in\mathcal P^{k+1}$.
\end{proposition}

\begin{proof} It follows from \eqref{ivds}
\[i(v)\mathcal DF=\lambda^1F^{-1}.j^1\pi_kF.v-v,\]
where $v\in\mathcal T$, so
\[\lambda^1F.(i(v)\mathcal DF)=j^1\pi_kF.v-\lambda^1F.v.\]
As $i(v)\mathcal DF\in \mathcal R^k$, we get $\lambda^1F.(i(v)\mathcal DF)\in \mathcal T\mathcal P^k$. Also, we get from $\pi_kF\in\mathcal P^k$ that $j^1\pi_kF.v\in\mathcal T\mathcal P^k$. Therefore,
\[\lambda^1F.v=j^1\pi_kF.v-\lambda^1F.(i(v)\mathcal DF)\in\mathcal T\mathcal P^k,\]
so $F\in\mathcal P^{k+1}$.
\end{proof}

If the linear Lie equation $R^k$ is formally integrable, and $P^k$ is the dif\/ferentiable groupoid associated with $R^k$, it is true (cf.\ Proposition~6.1, \cite{Ma2}) that the prolongation $P^{k+l}$ of $P^k$ is the groupoid associated with the linear Lie equation $R^{k+l}$. Therefore, $\pi_{k+l}:P^{k+l+1}\rightarrow P^{k+l}$ are submersions, for  $l\geq 0$.


\section{Formal isomorphism of intransitive linear Lie equations}\label{formaliso}
In the following sections, we consider intransitive linear Lie equations.
\begin{definition}
We say that a linear Lie equation $R^k\subset J^kT$ is \emph{intransitive} if there exists an integrable distribution $V\subset T$ such that $R^k\subset J^kV$ and $\pi_0(R^k)=J^0V$.
\end{definition}
In reality, considering \eqref{lietilde}, we need only to verify that $\pi_0(R^k)$ is a subvector bundle of~$J^0T$.
Our basic problem in this section is to determine the conditions for  two intransitive linear Lie equations to be isomorphic. This means  that there exists a dif\/feomorphism that sends one equation onto the other. In  the sequel, we give a brief description of the system of partial dif\/ferential equations that we should solve to obtain a class of dif\/feomorphisms $f:M\rightarrow M'$ such that $(j^{k+1}f)_*(R^k)=R'^k$. We utilize the same notation of
Section~\ref{partcon}. Consider $R^k\subset J^kV$ and $R'^k\subset J^kV'$ intransitive linear Lie equations, and $P^k\subset Q^k_V$ and $P'^k\subset Q^k_{V'}$ the associated groupoids.

\begin{definition}We say that a submanifold $S^k\subset Q^k_\phi$ is \emph{automorphic} by $P^k$ if $\alpha: S^k\rightarrow M$, $\beta:S^k\rightarrow M'$ are submersions, and for every  $X\in S^k(a,b)$, where $a\in M$ and $b\in M'$,
\[S^k(\cdot,b)=X\circ P^k(\cdot,a).\]
\end{definition}
We denote by $\mathcal S^k$ the set of invertible sections of $S^k$.
\begin{proposition} \label{pofor}Let $S^{k+1}$ be the prolongation of $S^k$. Then an invertible section $F\in \mathcal{Q}^{k+1}_\phi$ is such that $F(x)\in S^{k+1}(x)$ for every $x\in \alpha(F)$ if and only if
 $\pi_kF\in\mathcal{S}^k$ and $\mathcal{D}F\in\mathcal{T}^*\otimes\mathcal{R}^k.$
 \end{proposition}
 \begin{proof} The  same proof of  Proposition~\ref{prolqk} applies.
 \end{proof}

We def\/ine the \emph{symbol}
\[g^k_S=\{v\in TS^k:(\pi_{k-1})_*v=0\}.\]
The symbol $g^k_S$ of $S^k$ is isomorphic to the symbol $g^k$ of $R^k$, and we get an  complex analogous to~\eqref{deltas}, and we def\/ine that $g^k_S$ is 2-acyclic in the same way. From the formal integrability theorem (see \cite{G2}) we obtain:
\begin{proposition}
Let be  $S^k\subset Q^k_\phi$ automorphic by $P^k$  such that
\begin{enumerate}\itemsep=0pt
	\item[$(i)$] $S^{k+1}$ is a submanifold of $Q^{k+1}_\phi$,
	\item [$(ii)$]$\pi_{k}:S^{k+1}\rightarrow S^{k}$ is onto,
	\item [$(iii)$]$g^k_S$ is $2$-acyclic.
\end{enumerate}
Then $S^k$ is formally integrable, and each prolongation $S^{k+r}$ is automorphic by $P^{k+r}$, for $r\geq 1$.
\end{proposition}

\begin{definition}\label{foriso}
We say that the intransitive linear Lie equation $R^k\subset J^kV$ is \emph{formally isomorphic} to the intransitive linear Lie equation $R'^k\subset J^kV'$  at points $a$ and $a'$, respectively, if there exists a dif\/feomorphism
$\phi:N\rightarrow N'$, and
a submanifold $S^k\subset Q^k_\phi$ automorphic by $P^k$ and formally integrable, such that:
\begin{enumerate}\itemsep=0pt
	\item[$(i)$] $S'^k=\{X^{-1}:X\in S^k\}\subset Q^k_{\phi^{-1}}$ is automorphic by $P'^k$;
	\item[$(ii)$] $S^k(a,a')\neq \varnothing$.
\end{enumerate}
If there exists a solution $f:M\rightarrow M'$ of $S^k$, i.e., a dif\/feomorphism $f$ such that $j^kf$ is a section of $S^k$, and $f(a)=a'$, then $R^k$ at point $a$ is said \emph{isomorphic} to $R'^k$ at point $a'$.
\end{definition}
This def\/inition is essentially local. A most useful way to verify the formal isomorphism is given by proposition below, analogous of Proposition~\ref{pofor}:
 \begin{proposition}
 \label{FQ}
Suppose that $R^k\subset J^kV$, $R'^k\subset J^kV'$ are intransitive linear Lie equations, $N$~and~$N'$  submanifolds of $M$ and $M'$ transversal to integral submanifolds of $V$ and  $V'$, respectively, and $\phi:N\rightarrow N'$ a~diffeomorphism, $a\in N$, $a'\in N'$, and $\phi(a)=a'$.  Suppose furthermore that
the symbol $g^k$ of $R^k$ is $2$-acyclic. If there exists $F\in \mathcal{Q}^{k+1}_\phi$ such that $\beta F|_N=\phi$, $F_*(R^k)=R'^k$, and $\mathcal{D}F\in \mathcal{T}^*\otimes \mathcal{R}^k$, then $R^k$  at $a$ is formally isomorphic to $R'^k$ at $a'$.
\end{proposition}

\begin{proof} Def\/ine
\[ S^{k+1}=\big\{YF(x)X: X\in P^{k+1}(\cdot,x),\; Y\in P'^{k+1}(f(x),\cdot),\; x\in \alpha(F)\},\]
where $f=\beta F:\alpha(F)\subset M\rightarrow\beta(F)\subset M'$. Let be $U=\rho^{-1}(\rho(\alpha(F)))\subset M$ and $U'=\rho'^{-1}(\rho'(\beta(F)))\subset M'$. Observe that $S^{k+1}\subset Q^{k+1}_\phi$ and  $\alpha\times \beta:S^{k+1}\rightarrow U\times U'$ is onto (at least locally). If $S^k=\pi_kS^{k+1}$, then it is a  straightforward verif\/ication that $S^k$ is automorphic by $P^k$ and $S'^k$ is automorphic by $P'^k$.

Given an invertible section $G\in\mathcal{S}^{k+1}$, then in the neighborhood of each point of $\alpha(G)$, there are invertible sections $G_1\in\mathcal{P}^{k+1}$ and $G_2\in\mathcal{P'}^{k+1}$ such that $G=G_2FG_1$. In fact, given a point $x\in\alpha(G)$, there is an open set $V_x\subset \alpha(G)$, with $x\in V_x$, and an invertible section $G_1$ of $P^{k+1}$ def\/ined on $V_x$, such that $\beta(G_1)\subset \alpha(F)$. Let be $G_2=GG_1^{-1}F^{-1}$, def\/ined on $f(\beta(G_1))$. Then, $G_2$ is an invertible section of $P'^{k+1}$, and $G|_{V_x}=G_2FG_1$. It follows from Proposition~\ref{p44}$(i)$ and Proposition~\ref{pofor} that  $\mathcal{D}G\in \mathcal{T}^*\otimes \mathcal{R}^k$ on the open set $V_x$. As the $V_x$'s cover $\alpha(G)$, we get this property on all $\alpha(G)$. Therefore, $S^k$ is formally integrable, and conditions of Def\/inition~\ref{foriso} are satisf\/ied. \end{proof}

\begin{corollary}\label{fif}
Suppose that $R^k\subset J^kV$, $R'^k\subset J^kV'$ are intransitive linear Lie equations, $N$~and~$N'$  submanifolds of $M$ and $M'$, transversal to integral submanifolds of $V$ and  $V'$, respectively. Let be $\phi:N\rightarrow N'$ a diffeomorphism, $a\in N$, $a'\in N'$, and $\phi(a)=a'$.  Suppose furthermore that
the symbol $g^k$ of $R^k$ is $2$-acyclic. If there exists $F\in \mathcal{Q}^{k+1}_\phi$ such that $\beta F|_N=\phi$, and $F_*(T\oplus R^k)=T'\oplus R'^k$,  then $R^k$  at $a$ is formally isomorphic to $R'^k$ at $a'$.
\end{corollary}
Let's now show the existence of a f\/lat partial connection that leaves $R^k$ invariant.
\begin{proposition}\label{conrk}
Let be $R^k\subset J^kV$  an intransitive  
linear Lie equation. Then there exists a flat partial connection
\[\nabla: \ \mathcal{H}\oplus J^k\mathcal{V}\rightarrow \mathcal{V}^*\otimes(\mathcal{H}\oplus J^k\mathcal{V}),\]
such that, restricted to $\mathcal{H}\oplus \mathcal{R}^k$, it satisfies
\[\nabla: \ \mathcal{H}\oplus\mathcal{R}^k \rightarrow\mathcal{V}^*\otimes (\mathcal{H}\oplus\mathcal{R}^k).\]
Furthermore, if $U+\Xi$ is a  parallel section of $H\oplus J^kV$ and $(U+\Xi)|_N$ is a section of $(H\oplus R^k)|_N$, then $U+\Xi\in\mathcal{H}\oplus\mathcal{R}^k$.
\end{proposition}
\begin{proof} Choose a family of dif\/ferentiable sections $\sigma^y$ introduced in \eqref{sigmay} satisfying $\sigma^y(N)\subset P^{k+1}$. As $P^{k+1}$ is the groupoid associated with $R^{k+1}$, the form $\omega$ def\/ined by \eqref{omegacon} belongs to $\mathcal{V}^*\otimes\mathcal{R}^{k+1}$, and the partial connection $\nabla:\mathcal{H}\oplus J^k\mathcal{V}\rightarrow \mathcal{V}^*\otimes(\mathcal{H}\oplus J^k\mathcal{V})$
def\/ined by $\tilde\omega$, restricted to $H\oplus R^k$ sends $\mathcal{H}\oplus\mathcal{R}^k$ to $\mathcal{V}^*\otimes (\mathcal{H}\oplus\mathcal{R}^k)$, as a consequence of \eqref{con}. The proof now follows from Propositions~\ref{flat} and~\ref{UXI}.
\end{proof}

Now we prove the fundamental theorem for formal isomorphism of linear Lie equations:

\begin{theorem}\label{iso}
Suppose that $R^k\subset J^kV$ and $R'^k\subset J^kV'$ are intransitive linear Lie equations, $N$~and~$N'$   submanifolds of $M$ and $M'$ transversal to integral submanifolds of $V$ and  $V'$, respectively, and $\phi:N\rightarrow N'$ a diffeomorphism.  Suppose furthermore that  there exists  $\Phi:N\rightarrow \mathcal{Q}^{k+1}_\phi$ such that $\beta\Phi=\phi$, and
\[\Phi_*\big(TN\oplus R^k|_N\big)=TN'\oplus R'^k|_{N'}.\]
Then given a diffeomorphism $f:M\rightarrow M'$ such that $f_*V=V'$, $f|_N=\phi$, there exists $F\in\mathcal{Q}^{k+1}_\phi$ satisfying $F|_N=\Phi$, $\beta F=f$, and
\[F_*\big(T\oplus R^k\big)=T'\oplus R'^k.\]
\end{theorem}
\begin{proof} Let be  families $\sigma^y:N\rightarrow P^{k+1}$, $\sigma'^{y'}:N'\rightarrow P'^{k+1}$, as in the proof of Proposition~\ref{conrk}, and def\/ining f\/lat partial connections
\[\nabla: \ \mathcal{H}\oplus J^k\mathcal{V}\rightarrow \mathcal{V}^*\otimes \big(\mathcal{H}\oplus J^k\mathcal{V}\big),\]
 and
 \[\nabla': \ \mathcal{H}'\oplus J^k\mathcal{V}'\rightarrow \mathcal{V'}^*\otimes \big(\mathcal{H'}\oplus J^k\mathcal{V}'\big)\]
such that
\[\nabla\big(\mathcal{H}\oplus \mathcal{R}^k\big)\subset \mathcal{V}^*\otimes \big(\mathcal{H}\oplus \mathcal{R}^k\big)\]
and
\[\nabla'\big(\mathcal{H'}\oplus \mathcal{R'}^k\big)\subset \mathcal{V'}^*\otimes \big(\mathcal{H'}\oplus \mathcal{R'}^k\big).\]
 Observe that by Proposition~\ref {conrk}, $\omega\in\mathcal{V}^*\otimes \mathcal{R}^{k+1}$, and $\omega'\in\mathcal{V'}^*\otimes \mathcal{R'}^{k+1}$. Redef\/ine $H'=f_*H$, if necessary, to obtain $(x',y')=f(x,y)=(a(x),b(y))$, and def\/ine $F\in\mathcal{Q}^{k+1}_\phi$ by \[F(x,y)=\sigma'^{y'}(x',0)\Phi(x,0)\sigma^y(x,0)^{-1}.\]
  Then, we get from \eqref{sigmalinha}
\[F(x,y)\sigma^y(x,0)F(x,0)^{-1}=\sigma'^{y'}(x',0),\]
and from \eqref{nabla1} we get $\nabla'F_*=F_*\nabla$. By hypothesis
\[F_*\big(TN\oplus R^k|_N\big)=TN'\oplus R'^k|_{N'},\]
then by Proposition~\ref{conrk} we obtain
\begin{equation}\label{ftnrk}F_*\big(H\oplus R^k\big)=H'\oplus R'^{k}.
\end{equation}
From this and Corollary~\ref{sigmaD} we obtain
\begin{equation}\label{dhf}
\mathcal D_HF\in\mathcal H^*\otimes \mathcal R^k.
\end{equation}
It follows from \eqref{pik} that
\begin{equation*}
\pi_k\omega'=F_*(\mathcal{D}_VF+\pi_k\omega).
\end{equation*}
So, from $\pi_k\omega\in\mathcal{R}^k$, $\pi_k\omega'\in\mathcal{R'}^k$ and \eqref{ftnrk} we get
\[\mathcal{D}_VF\in\mathcal{V}^*\otimes \mathcal{R}^k.\]
Combining this with \eqref{ftnrk} and \eqref{dhf}, we get $\mathcal{D}F\in\mathcal{T}^*\otimes\mathcal{R}^k$ and $F_*(R^k)=R'^k$, and by Proposition~\ref{FQ} the theorem follows.
\end{proof}

\begin{corollary}\label{cofi}
Suppose that $R^k\subset J^kV$ and $R'^k\subset J^kV'$ are intransitive linear Lie equations, that $N$ and $N'$ are submanifolds of $M$ and $M'$ transversal to integral submanifolds of $V$ and  $V'$, respectively. Let be $\phi:N\rightarrow N'$ a diffeomorphism such that  $\phi(a)=a'$, where  $a\in N$ and $a'\in N'$.  Suppose furthermore that
the symbol $g^k$ of $R^k$ is $2$-acyclic.  If there exists $\Phi:N\rightarrow \mathcal{Q}^{k+1}_\phi$ such that $\beta\Phi=\phi$, and $\Phi_*(TN\oplus R^k|_N)=TN'\oplus R'^{k}|_{N'}$, then $R^k$ at point $a$ is formally isomorphic to $R'^k$ at point $a'$.
\end{corollary}
\begin{proof} The corollary follows from Theorem~\ref{iso} and Corollary~\ref{fif}.
\end{proof}

\section{Intransitive Lie algebras}\label{ila}
In this section, we associate an intransitive Lie algebra with a germ of an intransitive linear Lie equation.  This def\/inition must generalize the def\/inition of transitive Lie algebra, and incorporate the fact that we can reconstruct an intransitive linear Lie equation from its restriction to a~transversal to the orbits, unless of formal isomorphism, as the Theorem of \cite{V} and Theorem~\ref{iso} above shows.

We continue, in this section, to suppose
that $R^k\subset J^kV$ is an intransitive linear Lie equation and $g^k$ is $2$-acyclic. We remember that it follows from these hypotheses, see Corollary~\ref{rkfi},  that the prolongations $R^{k+l}$ of $R^k$, $l\geq 1$, satisfy:
\begin{enumerate}\itemsep=0pt
	\item[$(i)$] $R^{k+l}$ is a subvector bundle of $J^{k+l}V$;
	\item[$(ii)$] $\pi_{k+l}:R^{k+l+1}\rightarrow R^{k+l}$ is onto;
	\item[$(iii)$] $\firstbra{\mathcal T\oplus \mathcal R^{k+l}}{\mathcal T\oplus \mathcal R^{k+l}}_{k+l}
\subset{\mathcal T\oplus \mathcal R^{k+l-1}}$.
\end{enumerate}
We also assume  that $R^l=\pi_lR^k$ is a subvector bundle of $J^lV$ for every $0\leq l\leq k-1$, in particular, $R^0=J^0V$.

We denote by $\mathcal{O}_{N,a}$ the $\mathbb{R}$-algebra of germs at point $a\in N$ of local $C^\infty$ real functions on $N$. The $\mathcal{O}_{N,a}$-module $(\mathcal{T}N)_a$  of germs at point $a$ of local sections of $TN$ is isomorphic to $\mbox{Der }\mathcal{O}_{N,a}$, the $\mathcal{O}_{N,a}$-module of derivations of $\mathcal{O}_{N,a}$.
We denote by $L_j$ the $\mathcal{O}_{N,a}$-module of germs at point $a$ of local $C^\infty$-sections of $R^j|_N$, considered as a vector bundle on $N$ by the map $\pi|_N:R^k|_N\rightarrow N$. Let be the $\mathcal{O}_{N,a}$-modules
\[L=\lim\mbox{proj } L_j\]
and
\[\mathcal{L}=\mbox{Der }\mathcal{O}_{N,a}\oplus L.\]
The bilinear antisymmetric map
\[\firstbra{\;}{\;}_{k+l}: \ (\mathcal T\oplus \mathcal R^{k+l})\times(\mathcal T\oplus \mathcal R^{k+l})\rightarrow \mathcal T\oplus \mathcal R^{k+l-1}\]
induces a well def\/ined $\mathbb{R}$-bilinear antisymmetric map
\[\firstbra{\;}{\;}_j: \ (\mbox{Der }\mathcal{O}_{N,a}\oplus L_j )\times (\mbox{Der }\mathcal{O}_{N,a}\oplus L_j  )\rightarrow \mbox{Der }\mathcal{O}_{N,a}\oplus L_{j-1}.\]
As we saw in \eqref{checkcolch}, the projective limit of $\firstbra{\;}{\;}_j$ induces
a $\mathbb{R}$-Lie bracket
\[\firstbra{\;}{\;}_\infty: \ \mathcal{L}\times\mathcal{L}\rightarrow \mathcal{L},\]
so that $\mathcal{L}$ is a $\mathbb{R}$-Lie algebra.

The structure of the Lie algebra $\mathcal{L}$ is the semi-direct product of the Lie algebras $\mbox{Der }\mathcal{O}_{N,a}$ and~$L$, where the action of $\mbox{Der }\mathcal{O}_{N,a}$ on $L_m$ is given by the restriction  to $N$ of
$D_H:\mathcal{R}^j\rightarrow \mathcal{H}^*\otimes \mathcal{R}^{j-1}.$ If $v\in \mbox{Der }\mathcal{O}_{N,a}$ and $\xi$ is a section of $R^j|_N$ def\/ined in a neighborhood of $a\in N$, we get (see  Proposition~\ref{colchetek})
\[\firstbra{v}{\xi}_j=i(v)D\xi.\]
The map $(\rho_1)_*:\mathcal{L}\rightarrow \mbox{Der }\mathcal{O}_{N,a}$ is the canonical projection  given by the direct sum, and
\[\firstbra{\xi}{f\eta}_\infty=((\rho_1)_*\xi)(f)\eta+f\firstbra{\xi}{\eta}_\infty,\]
where $f\in\mathcal{O}_{N,a}$, $\xi,\eta\in \mathcal{L}$. The restriction of $\firstbra{\;}{\;}_\infty$ to $L$ is $\mathcal{O}_{N,a}$-bilinear, so $L$ is a $\mathcal{O}_{N,a}$-Lie algebra.  Each $L_j$ is a free $\mathcal{O}_{N,a}$-module f\/initely generated.

\begin{definition}
We call $\mathcal{L}$ the $\mathcal{O}_{N,a}$-\emph{intransitive $\mathbb{R}$-Lie algebra} associated with the formally integrable linear Lie equation $R^k$ at the point $a$ (and transversal $N$).
\end{definition}

In particular, we denote by $\mathcal{D}(V)$ the $\mathcal{O}_{N,a}$-intransitive Lie algebra associated with the linear Lie equation $J^0V$, and call it the $\mathcal{O}_{N,a}$-\emph{intransitive $\mathbb{R}$-Lie algebra associated with the involutive distribution $V$} at point $a\in M$. Clearly $\mathcal{L}\subset\mathcal{D}(V)$.

If $\mathcal{L}_j=\mbox{Der }\mathcal{O}_{N,a}\oplus L_j$, then $(\mathcal{L}_j,\firstbra{\;}{\;}_j)$ is called the \emph{truncated $\mathcal{O}_{N,a}$-intransitive $\mathbb{R}$-Lie algebra of order $j$} associated with $R^k$ at point $a$ (and transversal $N$). Then we can state the Theorem of \cite{V} as:
\begin{theorem}\label{VV}
Let be $\mathcal{L}_{k+2}\subset \mathcal{D}_{k+2}(V)$ a truncated $\mathcal{O}_{N,a}$-intransitive $\mathbb{R}$-Lie algebra. Then there exists a vector sub-bundle $R'^{k+1}\subset J^{k+1}V$ such that:
\begin{enumerate}\itemsep=0pt
	\item[$(i)$] $R'^k=\pi_k R'^{k+1}$ is a vector sub-bundle of $J^kV$;
	\item[$(ii)$] $\firstbra{\mathcal{T}\oplus \mathcal{R'}^{k+1}}{\mathcal{T}\oplus \mathcal{R'}^{k+1}}_{k+1}\subset \mathcal{T}\oplus \mathcal{R'}^{k}$;
	\item[$(iii)$] the truncated $\mathcal{O}_{N,a}$-intransitive  $\mathbb{R}$-Lie algebra associated with $R'^{k+1}$ is  $\mathcal{L}_{k+1}=\pi_{k+1}\mathcal{L}_{k+2}$.
\end{enumerate}
Furthermore, if $h_k=\{\xi\in\mathcal{L}_k|\pi_{k-1}\xi=0\}$ is $2$-acyclic, then $R'^k$ is formally integrable.
\end{theorem}

The def\/initions of $\mathcal{L}$ and $\mathcal{D}(V)$ depend on the choice  of the transversal $N$. Let's now introduce a notion of isomorphism inspired in Theorem~\ref{iso} such that  the intransitive Lie algebras obtained at point $a$ taking dif\/ferent transversal submanifolds are isomorphic.  We maintain the notation of Section~\ref{formaliso}.


Suppose that $R^k\subset J^kV$, $R'^k\subset J^kV'$ are formally integrable intransitive linear Lie equations, $N$ and $N'$  submanifolds of $M$ and $M'$ transversal to integral submanifolds of $V$ and  $V'$, respectively, and $\phi:N\rightarrow N'$ a dif\/feomorphism, $a\in N$, $a'\in N'$, and $\phi(a)=a'$.  We denote also by $\phi$ the isomorphism of $\mathbb{R}$-algebras $\phi:\mathcal{O}_{N,a}\rightarrow \mathcal{O}_{N',a'}$, def\/ined by $\phi(f)=f\phi^{-1}$. Let be
\[\Phi_{j+1}: \ N\rightarrow \mathcal{Q}^{j+1}_{\phi}\]
$\alpha$-sections
 such that $\phi=\beta\Phi_{j+1}$ and $\pi_j\Phi_{j+1}=\Phi_j$ for $j\geq 0$. Put $\Phi=\lim\mbox{proj }\Phi_j$.   We get maps
\[(\Phi_{j+1})_*: \ \mathcal{D}_j(V)\rightarrow \mathcal{D}_j(V')\]
and
\[\Phi_*: \ \mathcal{D}(V)\rightarrow \mathcal{D}(V').\]
The map $\Phi_*$ is $\mathbb{R}$-linear, commutes with $\firstbra{\;}{\;}_\infty$, and, if $f\in \mathcal{O}_{N,a}$ and $\xi\in\mathcal{D}(V)$, then \[\Phi_*(f\xi)=\phi(f)\Phi_*(\xi).\]
\begin{definition}
We say that $\Phi_*$ is an \emph{isomomorphism} from intransitive Lie algebra $\mathcal{L}\subset\mathcal{D}(V)$ \emph{onto} intransitive Lie algebra $\mathcal{L}'\subset\mathcal{D}(V')$ if
\[\Phi_*\mathcal{L}=\mathcal{L'}.\]
 If  $\Phi_*$ is an isomorphism, then $\mathcal{L}$ is said \emph{isomorphic} to $\mathcal{L}'$.
\end{definition}

\begin{proposition}Suppose that $R^k\subset J^kV$ is a formally integrable intransitive linear Lie equation, $a$, $b$  points of $M$, $N$, $N_1$ transversal  to the orbits of $R^k$ through the points $a$, $b$,
and~$\mathcal{L}$,~$\mathcal{L}_1$ the intransitive Lie algebras  associated with $R^k$ at the points $a$, $b$ $($and transversal $N$, $N_1)$, respectively. Let be $\rho:M\rightarrow N$ the fibration $($at least locally$)$ defined by the leaves of $V$. If $\rho(a)=\rho(b)$, then the $\mathcal{O}_{N,a}$-intransitive Lie algebra $\mathcal{L}$ is isomorphic to the $\mathcal{O}_{N_1,b}$-intransitive Lie algebra $\mathcal{L}_1$.
\end{proposition}

\begin{proof} If $x,y\in M$, $\rho(x)=\rho(y)$, there exists $X\in P^m$ with $\alpha(X)=x$, $\beta(X)=y$. Therefore, we can choose $\Phi_j:N\rightarrow P^j$ such that $\beta\Phi_j(N)=N_1$. Then
\[(\Phi_j)_*\big(TN\oplus R^{j-1}|_N\big)=TN_1\oplus R^{j-1}|_{N_1}.\]
It follows from the formal integrability of $P^k$ that we can choose the family $\{\Phi_j: j\geq 1\}$ such that $\Phi_{j+1}$ projects on $\Phi_j$, for $j\geq 1$. If $\Phi=\lim\mbox{proj }\Phi_j$, then $\Phi_*\mathcal{L}=\mathcal{L}_1$.
\end{proof}

With these def\/initions, we can state Corollary~\ref{cofi}  as:
\begin{theorem}\label{isore}Suppose that $R^k\subset J^kV$ is a linear Lie equation with symbol $g^k$ $2$-acyclic, and $R'^k\subset J^kV'$ another linear Lie equation. Let be  $a\in M$, $a'\in M'$, $N$ and $N'$ transversal to the orbits of $R^k$ and $R'^k$ through the points $a$ and $a'$, $\mathcal{L}_k$ and $\mathcal{L'}_k$  the truncated intransitive Lie algebras associated with $R^k$ and $R'^k$, at points $a$, $a'$ and tranversal $N$ and $N'$, respectively. If there exists $\Phi:N\rightarrow \mathcal{Q}^{k+1}_\phi$ such that $\beta\Phi=\phi:N\rightarrow N'$, $\phi(a)=a'$, and $\Phi_*\mathcal{L}_k=\mathcal{L'}_k$, then $R^k$ at point $a$ is formally isomorphic to $R'^k$ at point $a'$.
\end{theorem}

\section{Application}
As an application of this theory, we could utilize the def\/inition of intransitive Lie algebras to obtain the intransitive linear Lie equations in the plane obtained by \'E.~Cartan in \cite{C}. We will be limited to  classifying the f\/irst order intransitive linear Lie equations, with $\dim g^1=1$. This will include  the example we presented in the introduction, which was not presented by Cartan in his table, suppressed by a nullity hypothesis.

Let be $V$  a 1-dimensional distribution on $\mathbb{R}^2$, which we can suppose is generated by the vector f\/ield $\frac{\partial}{\partial y}$. We will use  the coordinate system $p_{j,l}$, $j,l\geq 0$, $0\leq j+l\leq k$, in $J^kV$, def\/ined by
\[p_{j,l}(j^k_{(a,b)}\Theta)=\frac{\partial^{j+l}\theta}{\partial x^j\partial y^l}{(a,b)},\]
where $\Theta(x,y)=\theta(x,y)\frac{\partial}{\partial y}$.
Let's consider $(0,0)$ as point base  and the transversal $N=\mathbb{R}\times\{0\}$. Then $\mbox{Der }\mathcal{O}_{N,(0,0)}$ is generated, as $\mathcal{O}_{N,(0,0)}$-module, by $\frac{\partial}{\partial x}|_N$.
Let be $g^k_V$ the symbol of $J^kV$. This symbol is generated by $f^{j,l}\otimes j^0\frac{\partial}{\partial y}$, where
\[f^{j,l}=\frac{1}{j!l!}(dx)^j(dy)^l,\]
with $j,l\geq 0$ and $j+l=k$ (cf.~\eqref{f}).
 Then
\[\left[\!\!\left[\frac{\partial}{\partial x},f^{j,l}\otimes j^0\frac{\partial}{\partial y}\right]\!\!\right]_{k}=-f^{j-1,l}\otimes j^0\frac{\partial}{\partial y}\in g^{k-1}_V,\]
and if $Y\in J^kV$ is such that $\pi_0(Y)=j^0\frac{\partial}{\partial y}$, then
\[\left[\!\!\left[Y,f^{j,l}\otimes j^0\frac{\partial}{\partial y}\right]\!\!\right]_k=f^{j,l-1}\otimes j^0\frac{\partial}{\partial y}\in g^{k-1}_V.\]
Let be $R^1\subset J^1V$ a formally integrable linear Lie equation, with $R^0=J^0V$. Suppose $g^1$ is the symbol of $R^1$, with $\dim g^1=1$. Then $g^1$ is generated by an element $X_1=(Af^{1,0}+Bf^{0,1})\otimes j^0\frac{\partial}{\partial y}$, with $A(x,y)^2+B(x,y)^2\neq 0$.
 \begin{lemma}
 If $R^2$ is the prolongation of $R^1$, then the symbol $g^2$ of $R^2$ is generated by
 \[\big(Af^{1,0}+Bf^{0,1}\big)^2\otimes j^0\frac{\partial}{\partial y}.\]
\end{lemma}

\begin{proof} Consider $Y_0,Y_1\in \mathcal R^2$ such that $\pi_0Y_0=j^0\frac{\partial}{\partial y}$ and $\pi_1Y_1=(Af^{1,0}+Bf^{0,1})\otimes j^0\frac{\partial}{\partial y}$ is a~section of $g^1$. Let be $Y=(a_{20}f^{2,0}+a_{11}f^{1,1}+a_{02}f^{0,2})\otimes j^0\frac{\partial}{\partial y}$ a section of $g^2$. Then
\[\left[\!\!\left[\frac{\partial}{\partial x},Y\right]\!\!\right]_{2}
=-\big(a_{20}f^{1,0}+a_{11}f^{0,1}\big)\in g^{1},\]
so $(a_{20}f^{1,0}+a_{11}f^{0,1})=\lambda(Af^{1,0}+Bf^{0,1})$, for  a real function $\lambda$. In a similar way,
\[\firstbra{Y_0}{Y}_{2}=\big(a_{11}f^{1,0}+a_{02}f^{0,1}\big)\in g^{1},\]
and $(a_{11}f^{1,0}+a_{02}f^{0,1})=\mu(Af^{1,0}+Bf^{0,1})$, for some real function $\mu$. Therefore, we get $a_{20}=\lambda A$, $a_{11}=\lambda B=\mu A$ and $a_{02}=\mu B$. So, there exists $r$ such that $\lambda=rA$ and $\mu=rB$. From this, $a_{20}=r A^2$, $a_{11}=rAB$ and $a_{02}=rB^2$. Then
\begin{gather*}
Y=r\big(Af^{1,0}+Bf^{0,1}\big)^2\otimes j^0\frac{\partial}{\partial y}.\tag*{\qed}
\end{gather*}
\renewcommand{\qed}{}
\end{proof}

A similar argument shows that $g^k$ is one dimensional and is generated by $(Af^{1,0}+Bf^{0,1})^k\otimes j^0\frac{\partial}{\partial y}.$ Consider now the complex \eqref{deltas},
\[
0\rightarrow g^{l+1}\stackrel{\delta}{\rightarrow} T^*\otimes g^{l}\stackrel{\delta}{\rightarrow}\wedge^2\mathcal T^*\otimes g^{l-1} {\rightarrow}0
\]
for $l\geq 2$. As $\dim g^k=1$ and $\dim T=2$, this complex is clearly exact, so $g^1$ is 2-acyclic.

Let's now verify conditions on the  truncated intransitive Lie algebra $\mathcal L_2$. The $\mathcal{O}_{N,(0,0)}$-module  $\mathcal TN_{(0,0)}$ is generated by \[Y_{-1}=\frac{\partial}{\partial x}\Big|_N,\]
and the generators of $L_2$ are $Y_0$, $Y_1$, def\/ined in the proof of lemma and $Y_2=\frac 1 2(af^{1,0}+bf^{0,1})^2\otimes j^0\frac{\partial}{\partial y}$, restricted to $N$,  again denoted by the same letters. Here,  $a(x)=A(x,0)$ and $b(x)=B(x,0)$. We have
\begin{alignat*}{3}
& \firstbra{Y_{-1}}{Y_0}_2=b_{00}\pi_1(Y_0)+b_{01}\pi_1(Y_1),\qquad && \firstbra{Y_{0}}{Y_1}_2=b\pi_1(Y_0)+a_{01}\pi_1(Y_1),&\\
& \firstbra{Y_{-1}}{Y_1}_2=-a\pi_1(Y_0)+b_{11}\pi_1(Y_1),\qquad && \firstbra{Y_{0}}{Y_2}_2=b\pi_1(Y_1),& \\
& \firstbra{Y_{-1}}{Y_2}_2=-a\pi_1(Y_1),\qquad && \firstbra{Y_{1}}{Y_2}_2=0,&
\end{alignat*}
for $a$, $b$, $b_{00}$, $b_{01}$, $b_{11}$, $a_{01}$  in $\mathcal{O}_{N,(0,0)}$. It follows from
\[
\firstbra{Y_{-1}}{\firstbra{Y_0}{Y_1}_2}_1=\firstbra{\firstbra{Y_{-1}}{Y_0}_2}{\pi_1(Y_{1})}_1+\firstbra{\pi_1(Y_0)}{\firstbra{Y_{-1}}{Y_{1}}_2}_1
\]
that $(\frac{\partial b}{\partial x}-aa_{01}-bb_{11})\pi_0(Y_0)=0$, so
\[\frac{\partial b}{\partial x}-aa_{01}-bb_{11}=0.\]
As $a(0)$ or $b(0)$ is not null, we can solve this equation for $a_{01}$ or $b_{11}$. Then we can  f\/ind a~truncated intransitive Lie algebra $\mathcal L_2$ that projects on $\mathcal L_1$. Now by Theorem~\ref{isore} we classify the isomorphism class of $\mathcal L_1$. We must examine two cases:

\textbf{Case 1.} We suppose $b(0)\neq 0$. By dividing by $b$, we can suppose $b=1$. Then the algeb\-ra~$\mathcal L_1$ is generated by $Y_{-1}$, $\pi_0Y_0$, $\pi_1Y_1$, and they satisfy, from above,
\begin{gather*}
\firstbra{Y_{-1}}{\pi_1Y_0}_1=b_{00}\pi_0(Y_0),\qquad
\firstbra{Y_{-1}}{\pi_1Y_1}_1=-a\pi_0(Y_0),\qquad
\firstbra{\pi_1Y_{0}}{\pi_1Y_1}_1=\pi_0(Y_0).
\end{gather*}
Let be $\mathcal L'_1$ a truncated intransitive Lie algebra generated by $X_{-1}$, $X_0$ and $X_1$, such that
\begin{gather*}
\firstbra{X_{-1}}{X_0}_1=0,\qquad
\firstbra{X_{-1}}{X_1}_1=0,\qquad
\firstbra{X_{0}}{X_1}_1=\pi_0(X_0),
\end{gather*}
and $f_1:\mathcal L'_1\rightarrow \mathcal L_1$ def\/ined by
\begin{gather*}
f_1(X_{-1})=Y_{-1}+a\pi_1(Y_0)+b_{00}\pi_1(Y_1),\qquad
f_1(X_0)=\pi_1(Y_0),\qquad
f_1(X_1)=\pi_1(Y_1).
\end{gather*}
Let's verify that $f_1$ is an isomorphism of truncated intransitive Lie algebras of order 1. In fact,
\begin{gather*}
\firstbra{f_1(X_{-1})}{f_1(X_{0})}_1=\firstbra{Y_{-1}+a\pi_1(Y_0)+b_{00}\pi_1(Y_1)}{\pi_1(Y_0)}_1=0
=f_0(\firstbra{X_{-1}}{X_0}_1),\\
 \firstbra{f_1(X_{-1})}{f_1(X_{1})}_1=\firstbra{Y_{-1}+a\pi_1(Y_0)+b_{00}\pi_1(Y_1)}{\pi_1(Y_1)}_1=0
 =f_0(\firstbra{X_{-1}}{X_1}_1),
 \end{gather*}
and
\[\firstbra{f_1(X_{0})}{f_1(X_{1})}_1=\firstbra{\pi_1(Y_0)}{\pi_1(Y_1)}_1=\pi_0(Y_0)=f_0( \pi_0(X_0))= f_0(\firstbra{X_{0}}{X_1}_1).\]
 The truncated Lie algebra $\mathcal L'_1$ generated as $\mathcal O_{N,(0,0)}$-module by $X_{-1}$, $X_0$, $X_1$ can be represented~by
 \[X_{-1}=\frac{\partial}{\partial x}\Big|_N, \qquad X_0=j^1\frac{\partial}{\partial y},\qquad X_1=f^{0,1}\otimes j^0\frac{\partial}{\partial y}.\]
 The linear Lie equation associated with $\mathcal L'_1$ is
\[R'^1=\big\{(p_{0,0},p_{1,0},p_{0,1})\in J^1V: p_{1,0}=0 \big\},\]
so $R^1$ is formally isomorphic to $R'^1$.
The inf\/initesimal pseudogroup of solutions of $R'^1$ is
\[\left\{\Theta(x,y)=\theta(y)\frac{\partial}{\partial y}\right\}.\]

\textbf{Case 2.} We suppose $b(0)= 0$. In this case, $a(0)\neq 0$, and, dividing it by $a$, we can suppose $a=1$. Then the truncated intransitive Lie algebra $\mathcal L_1$ is generated by $Y_{-1}$, $\pi_0Y_0$, $\pi_1Y_1$, and they satisfy, from above,
\begin{gather*}
\firstbra{Y_{-1}}{\pi_1Y_0}_1=b_{00}\pi_0(Y_0),\qquad
\firstbra{Y_{-1}}{\pi_1Y_1}_1=-\pi_0(Y_0),\qquad
\firstbra{\pi_1Y_{0}}{\pi_1Y_1}_1=b\pi_0(Y_0).
\end{gather*}
Replacing $\pi_1Y_0$ by $\pi_1Y_0+b_{00}\pi_1Y_1$, we obtain $\firstbra{Y_{-1}}{\pi_1Y_0+b_{00}\pi_1Y_1}_1=0$, and the other products remain unchanged. Without loss of generality, we can suppose $b_{00}=0$. Let be $\mathcal L'_1$ a truncated intransitive Lie algebra generated by $X_{-1}$, $X_0$ and $X_1$, such that
\begin{gather*}
\firstbra{X_{-1}}{X_0}_1=0,\qquad
\firstbra{X_{-1}}{X_1}_1=-\pi_0(X_0),\qquad
\firstbra{X_{0}}{X_1}_1=\beta\pi_0(X_0),
\end{gather*}
and $f_1:\mathcal L'_1\rightarrow \mathcal L_1$ def\/ined by
\begin{gather*}
f_1(X_{-1})=c^{-1}Y_{-1},\qquad
f_1(X_0)=\pi_1(Y_0),\qquad
f_1(X_1)=c\pi_1(Y_1),
\end{gather*}
with $c(0)\neq 0$. Then
\begin{gather*}
\firstbra{f_1(X_{-1})}{f_1(X_{0})}_1=\firstbra{c^{-1}Y_{-1}}{\pi_1(Y_0)}_1=0=f_0(\firstbra{X_{-1}}{X_0}_1),\\
\firstbra{f_1(X_{-1})}{f_1(X_{1})}_1=\firstbra{c^{-1}Y_{-1}}{c\pi_1(Y_1)}_1=-\pi_0(Y_0)=f_0(\firstbra{X_{-1}}{X_1}_1),
\end{gather*}
and
\[\firstbra{f_1(X_{0})}{f_1(X_{1})}_1=\firstbra{\pi_1(Y_0)}{c\pi_1(Y_1)}_1=cb\pi_0(Y_0)=f_0(\beta \pi_0(X_0))= f_0(\firstbra{X_{0}}{X_1}_1),\]
if $cb=\beta.$

Then  $f_1$ is an isomorphism of truncated intransitive Lie algebras of order 1 if and only if there exists $c\in\mathcal O_{N,(0,0)}$ with $c(0)\neq 0$ such that $cb=\beta$.
The class of truncated intransitive Lie algebras is the class of equivalence of $b$, $b(0)=0$, and $b\equiv\beta$ if and only if there exists $c$, $c(0)\neq 0$, with $bc=\beta$.

 The truncated Lie algebra $\mathcal L'_1$ generated as $\mathcal O_{\mathbb{R}}$-module by $X_{-1}$, $X_0$, $X_1$ can be represented by
 \[X_{-1}=\frac{\partial}{\partial x},\qquad X_0=j^1\frac{\partial}{\partial y},\qquad X_1=\big(f^{1,0}+\beta f^{0,1}\big)\otimes j^0\frac{\partial}{\partial y},\]
 and the linear Lie equation associated with $\mathcal L'_1$ is
\[R'^1=\big\{(p_{0,0},p_{1,0},p_{0,1})\in J^1V: p_{0,1}=\beta p_{1,0} \big\},\]
so $R^1$ is formally isomorphic to $R'^1$.
The inf\/initesimal pseudogroup of solutions of $R'^1$ is
\[\left\{\Theta(x,y)=\theta(x,y)\frac{\partial}{\partial y}:\frac{\partial\theta}{\partial y}=\beta\frac{\partial\theta}{\partial x}\right\}.\]
If $\beta(x)=x$, we obtain
\[\left\{\Theta(x,y)=\theta(xe^y)\frac{\partial}{\partial y}\right\};\]
if $\beta(x)=x^k$, $k\geq 2$, we obtain
\[\left\{\Theta(x,y)=\theta\left(\frac{x^{k-1}}{(k-1)yx^{k-1}-1}\right)\frac{\partial}{\partial y}\right\};\]
if $\beta(x)=0$, we obtain
\[\left\{\Theta(x,y)=\theta(x)\frac{\partial}{\partial y}\right\}.\]
In the classif\/ication of \cite{C},  case 2 is represented only by  $\beta=0$.

\section{Conclusion}

The results of this paper show that the intransitive Lie algebra here introduced to represent a~linear Lie equation at a point is suf\/f\/icient to guarantee the existence and formal isomorphism of intransitive linear Lie equations. This
brings a new way to pursue the study of intransitive Lie groups and the applications envisaged by Sophus Lie on the integrability of partial dif\/ferential equations with a pseudogroup of invariants. It is clear that several problems can still exist,  as the relationship between   subalgebras of transitive algebras and intransitive algebras, and the notion of equivalence of intransitive algebras. A very interesting problem is the classif\/ication of simple intransitive Lie groups, since Cartan, in his list, excluded some classes of simple intransitive Lie group, as the example presented above.

\subsection*{Acknowledgments}

I would like to thank  the referees for the several suggestions to improve this paper.

\pdfbookmark[1]{References}{ref}
\LastPageEnding

\end{document}